\documentclass{article}

\usepackage{graphicx}%
\usepackage{amsmath,amssymb,amsfonts}%
\usepackage{amsthm}%
\usepackage{mathrsfs}%
\usepackage[title]{appendix}%
\usepackage{xcolor}%
\usepackage{tikz}

\theoremstyle{plain}
\newtheorem{prp}{Proposition}
\newtheorem{lem}{Lemma}

\newtheorem{cor}{Corollary}

\theoremstyle{remark}
\newtheorem{exa}{Example}
\newtheorem*{rem}{Remark}

\theoremstyle{definition}

\begin{document}
\title{Signed Minkowski decompositions of convex polygons into minimum simplices and factorization of max-plus functions}
\author{Soujun Kitagawa \\
skitagawa@aoni.waseda.jp \\
{Department of Fundamental Science and Engineering}, \\
{Waseda University}, \\
{1-104 Totsukamachi, Shinjuku-ku, {Tokyo}, {Japan}}}
\date{}

\maketitle
\abstract{Signed Minkowski decomposition is an expression of a polytope as a Minkowski sum and difference of smaller polytopes. Signed Minkowski decompositions of a polytope can be interpreted as factorizations of a max-plus (tropical) function. We review two relations about Minkowski decompositions, and we prove that any 2-dimensional integral polytopes (polygons) have a signed Minkowski decomposition which consists of integral points, integral line segments of length 1, and integral triangles of area 1/2. From this result, we also obtain that any max-plus functions with two variables and integer coefficients can be expressed as a set of specified form of simpler max-plus functions.}
\section{Introduction}
This paper concerns the decomposition of max-plus functions depending on two variables with integer coefficients. The decomposition of max-plus functions is expressed as the sum of generating max-plus functions, and is related to the Minkowski decomposition of convex polygons. This has been studied from various perspectives, such as computational geometry and algebraic combinatorics, and numerous results have been obtained~\cite{ABD2010, ASV2021,CFF2017,ET2006,Fukuda2004,Funke2021,GL2001,Hachenberger2009,JR2022,Lange2013}. In this paper, we prove that a convex polygon whose vertices are all at integral points always has a Minkowski decomposition composed of line segments of length 1 and integral triangles of area 1/2, where such triangles contain no integral points in their interior. A concrete procedure for obtaining such a decomposition is also derived. In particular, we will show that this decomposition can be expressed in terms of the sum and difference of elementary max-plus functions with a smaller number of terms. Such a decomposition of the max-plus function provides an effective expression for the improvement of neural networks~\cite{MRZ2022, TW2024} and for analyzing of time evolution equations in max-plus form~\cite{Nakata2021}.

For later reference, we now review the fundamental concepts of convex polytopes, Minkowski decompositions, and max-plus functions.

\textbf{Convex Polytopes}
Convex polytopes are fundamental objects in geometry, and arise in various branches of mathematics, including algebraic geometry~\cite{GKZ1994}, computational geometry~\cite{PS1985}, and combinatorics~\cite{Gruenbaum2003,Ziegler1995}. There are several classes of polytopes which are defined in terms of the information of some special structure, such as poset polytopes~\cite{Stanley1986} and generalized permutohedra~\cite{Postnikov2009}. In this paper, we consider Newton polytopes, which are defined by exponent vectors of polynomials, and reflect the algebraic properties of the associated polynomials. In fact, some classes of polynomials on the tropical semifield and their Newton polytopes are so closely related that they can be considered in parallel. The addition and product of tropical polynomials can be computed by taking convex hulls and Minkowski sums of the associated polytopes. Therefore, some results on polytopes can be directly rephrased as results on tropical polynomials. We consider problems for signed Minkowski decompositions of integral polytopes. Integral polytopes are polytopes whose vertices are on $\mathbb{Z}^n$, and the signed Minkowski decompositions, which were introduced by Ardila, Benedetti, and Doker~\cite{ABD2010}, are decompositions of polytopes with Minkowski sums and Minkowski differences.

\textbf{Minkowski decompositions}
The Minkowski sum for the subsets of $\mathbb{R}^2$ is obtained by the sum of position vectors for all combinations of elements. Expressing a convex polytope in terms of Minkowski sums of other convex polytopes is known as the Minkowski decomposition. The Minkowski sum and decomposition are actively studied in mathematical morphology~\cite{Serra1983}, since they are effective in expressing operations on numerical geometric data. A substantial amount of research has been undertaken into Minkowski decompositions from the prespective of applications, such as efficient algorithms and numerical approximation methods~\cite{ET2006,Fukuda2004,Hachenberger2009}. Minkowski decompositions are also of interest because of their mathematical structure. For example, it is known that if the polynomial is factorized, then its Newton polytope can be decomposed into the Newton polytopes of its factors~\cite{GL2001, Ostrowski1975, Ostrowski1976}.

The Minkowski difference is the inverse object of the Minkowski sum. The Minkowski difference does not always exist for any pair of polytopes. Thus, polytopes and their Minkowski sums do not form a group, but a monoid. The group completion (the Grothendieck group) of the monoid of polytopes is known as the group of virtual polytopes~\cite{PK1993,PS2015}. The subgroup of virtual polytopes consisting of integral polytopes have been studied by Cha, Friedl, Funke~\cite{CFF2017}, and Funke~\cite{Funke2021}. This subgroup is referred to as the integral polytope group.

The signed Minkowski decompositions of a polytope are defined as an expression of the polytope with Minkowski sums and differences. They are adopted in the study of matroid polytopes and generalized permutohedra~\cite{ABD2010,JR2022}. The partition of a polytope (also known as the polytopal subdivision) is a collection of smaller polytopes which are obtained by dividing the original polytope~\cite{Khovanskii1997, Ziegler1995}. Funke proved that any partition of a polytope induces an associated signed Minkowski decomposition~\cite{Funke2021}.

\textbf{Max-plus functions}
The union of real numbers $\mathbb{R}$ and $\{-\infty\}$, equipped with maximum function as addition and the usual summation as product, is known as max-plus algebra or tropical semifield. The max-plus algebra has been applied in various research fields, such as system analysis~\cite{BCOQ1992}, phylogenetics~\cite{BHV2001}, and neural networks~\cite{MCT2021,Robinson2021, ZNL2018}. The polynomials in a tropical semifield are called tropical polynomials. Tropical polynomials with $n$ variables as functions of $\mathbb{R}^n \rightarrow \mathbb{R}$ are continuous, piecewise linear, and convex. The study of tropical polynomials using tools of geometry, especially those of algebraic geometry, is known as tropical geometry~\cite{IMS2009,MS2015}. Newton polytopes of tropical polynomials are defined in the same manner as those of usual polynomials. Minkowski sums and differences of Newton polytopes correspond to products and divisions of the tropical polynomials. Thus, factorizations of tropical rational functions are induced from decompositions of polytopes.

The remainder of this paper is organized as follows: In Section~\ref{Sec:PolytopesAndConvexGeometry}, we review some terminology and definitions about convex polygons. In Section~\ref{Sec:PartitionsOfPolytopes}, we introduce identities about the Minkowski decomposition of convex polygons, known as `Cutting relations'. We then introduce partitions of convex polygons, and give `Partition relations' which are induced from a partition of a convex polygon by applying the cutting relations. In Section~\ref{Sec:DecompOfLatPolytopes}, we consider decompositions of integral polygons. We prove that any integral polygon can be expressed as a signed Minkowski decomposition consisting of unit segments, which are line segments of length 1, and unit triangles, which are integral triangles whose area is 1/2 (Proposition~\ref{PropOfDecompFormula}). In Section~\ref{Sec:maxplusfunctions}, we review max-plus functions, and examine the relationship between max-plus functions and geometry of convex polygons. In Section~\ref{Sec:Concluding}, we provide concluding remarks and further discussion.
\section{Polygons and convex geometry} \label{Sec:PolytopesAndConvexGeometry}
We begin by recalling some terminology and definitions about convex geometry and polygons. It is assumed that all of the objects under consideration are in Euclidean space $\mathbb{R}^2$. Most of the concepts in this paper can be applied in higher dimensional spaces, except Proposition~\ref{PropOfDecompFormula} in section~\ref{Sec:DecompOfLatPolytopes}. In this section, we introduce basic results without proofs. For complete proofs, see, e.g.,~\cite{DRS2010,Gruenbaum2003, Ziegler1995}. Two equivalent definitions of convex polytopes are known. One is defined as the convex hull of vertices, and the other is defined as the region enclosed by hyperplanes. We adopt the former one. In this paper, we often write polygon instead of convex polygon, since non-convex polygons do not appear.

Assume $P = \{\mathbf{p}_1 = (a_1, b_1), \dots, \mathbf{p}_m = (a_m, b_m)\}$ is a finite set of points in $\mathbb{R}^2$. Here $\mathbf{p}_i = (a_i,b_i) \in \mathbb{R}^2$ denotes both the coordinates of a point and its position vector. A subset $K$ of $\mathbb{R}^2$ is convex if for any two points $\mathbf{p}, \mathbf{p}' \in K$, the segment $\{t \mathbf{p} + (1 - t) \mathbf{p'} \,|\, 0 \le t \le 1\}$ between $\mathbf{p}$ and $\mathbf{p}'$ is contained in $K$. For a subset $S$ in $\mathbb{R}^2$, the convex hull of $S$, written as $\mathrm{conv}(S)$, is the minimal convex set including $S$. For a finite set $P$, we write $\mathrm{conv}(P) = \mathrm{conv}\{\mathbf{p}_1,\dots \mathbf{p}_m\}$ instead of $\mathrm{conv}(\{\mathbf{p}_1,\dots \mathbf{p}_m\})$. If a finite set $P$ is given, we can determine a unique subset $v(P) := \{\mathbf{p}_{v_1}, \dots, \mathbf{p}_{v_k}\}$ of $P$ such that $\mathrm{conv}(v(P)) = \mathrm{conv}(P)$ and $|v(P)| \le |P'|$ for any $P'$ satisfying $\mathrm{conv}(P') = \mathrm{conv}(P)$. The elements of $v(P)$ are called the vertices of $\mathrm{conv}(P)$. It is known that $\mathrm{conv}(P_1) = \mathrm{conv}(P_2)$ if and only if $v(P_1) = v(P_2)$. If $|v(P)| = 1$, $|v(P)| = 2$, or $|v(P)| \ge 3$, $\mathrm{conv}(P)$ is a point, a segment, or a polygon, respectively. Note that we will also use $v(\mathrm{conv}(P))$ to denote the vertices of $\mathrm{conv}(P)$. Additionally, $\mathrm{conv(P)}$ can be represented as
\begin{equation} \label{EquiOfVertex}
       \mathrm{conv}(P) = \{t_1 \mathbf{p}_{v_1} + \dots + t_k \mathbf{p}_{v_k} \,|\, t_i \ge 0,\, t_1 + \dots + t_k = 1\},
\end{equation}
with vertices $\mathbf{p}_{v_i} \in v(P)$.

A Minkowski sum $S_1 + S_2$ of two sets $S_1 \subset \mathbb{R}^2$ and $S_2 \subset \mathbb{R}^2$ is defined as
\begin{equation}
       S_1 + S_2 := \{\mathbf{s}_1 + \mathbf{s}_2 \,|\, \mathbf{s}_1 \in S_1, \mathbf{s}_2 \in S_2\},
\end{equation}
where $\mathbf{s}_1 + \mathbf{s}_2$ is a point obtained from the sum of the position vectors of $\mathbf{s}_1$ and $\mathbf{s}_2$. The Minkowski difference $S_1 - S_2$ between $S_1$ and $S_2$ can be defined in terms of the set $S_3$, if it satisfies $S_3 + S_2 = S_1$. Note that $S_1 - S_2$ is not the set $\{\mathbf{s}_1 + (-\mathbf{s}_2) \,|\, \mathbf{s}_1 \in S_1, \mathbf{s}_2 \in S_2\}$, and Minkowski differences cannot always be defined for every pair of $S_1$ and $S_2$. A Minkowski sum $S + \{\mathbf{p}\}$ between a set $S$ and a singleton $\{\mathbf{p}\}$ is obtained as the translation of $S$ by $\mathbf{p}$. We write $S + \mathbf{p}$ instead of $S + \{\mathbf{p}\}$. Similarly, we write $S - \mathbf{p}$ instead of $S - \{\mathbf{p}\}$, which is equal to $S + (-\mathbf{p})$.

The following lemma lists some useful identities, which are used throughout this paper. They can be obtained directly from the definitions.
\begin{lem} \label{LemOfGeometricIdentities}
       Assume $P_1$, $P_2$, and $P_3$ are finite sets of points on $\mathbb{R}^2$. Then, the following identities hold.
       \begin{align}
              (P_1 \cup P_2) + P_3 &= (P_1 + P_3) \cup (P_2 + P_3), \label{rel1}\\ 
              \mathrm{conv}(P_1 \cup P_2) &= \mathrm{conv}(\mathrm{conv} (P_1) \cup \mathrm{conv} (P_2)), \label{rel2}\\
              \mathrm{conv}(P_1 + P_2) &= \mathrm{conv}(P_1) + \mathrm{conv}(P_2), \label{rel3} \\
              P_1 - (P_2 + P_3) &= (P_1 - P_2) - P_3 = (P_1 - P_3) - P_2, \label{rel4}
       \end{align}
       where we assume that all of the differences appearing in \eqref{rel4} exist.
\end{lem}
We write $P_1 - P_2 - P_3 := P_1 - (P_2 + P_3) = (P_1 - P_2) - P_3 = (P_1 - P_3) - P_2$, since \eqref{rel4} holds. 
\section{Cutting relations and Partitions of polygons} \label{Sec:PartitionsOfPolytopes}
\subsection{Cutting relations of segments and polygons}
In this section, we will review two types of identities, known as the cutting relations~\cite{Funke2021}. Although these identities have been proven for higher dimensional polytopes, we review the results for segments and polygons only, which are sufficient for our purposes. For convenience, we use the notation such that $\overline{\mathbf{p}_{1}\mathbf{p}_{2}} := \mathrm{conv} \{\mathbf{p}_{1}, \mathbf{p}_{2}\}$ for $\mathbf{p}_1 \ne \mathbf{p}_2$. $\overline{\mathbf{p}_{1}\mathbf{p}_{2}}$ is a line segment with endpoints $\mathbf{p}_1$ and $\mathbf{p}_2$.
\begin{lem} \label{LemOfLine}
       Let $\mathbf{p}_1, \mathbf{p}_2 \in \mathbb{R}^2$, and $\mathbf{p}_3 := t \mathbf{p}_1 + (1-t) \mathbf{p}_2$, where $0 < t < 1$. Then the following identity holds.
       \begin{equation}
       \begin{gathered}
              \overline{\mathbf{p}_1 \mathbf{p}_2} =\overline{\mathbf{p}_1 \mathbf{p}_3} + \overline{\mathbf{p}_3 \mathbf{p}_2} - \mathbf{p}_3. \\
              \Leftrightarrow
              \begin{tikzpicture}[baseline = 1em, x = 1.4em, y = 1.4em]
                     \coordinate (A) at (0,0);
                     \coordinate (B) at (5,2);
                     \coordinate (C) at (3,1.2);
                     \draw[black] (A) -- (B);
                     \fill[black] (A) circle (0.11);
                     \fill[black] (B) circle (0.11);
                     \fill[black] (C) circle (0.11);
                     \node[above] at (A) {$\mathbf{p}_1$};
                     \node[below] at (B) {$\mathbf{p}_2$};
                     \node[below] at (C) {$\mathbf{p}_3$};
              \end{tikzpicture} =
              \begin{tikzpicture}[baseline = 1em, x = 1.4em, y = 1.4em]
                     \coordinate (A) at (0,0);
                     \coordinate (B) at (3,1.2);
                     \coordinate (C) at (5,2);
                     \draw[black] (A) -- (B);
                     \draw[dashed] (B) -- (C);
                     \fill[black] (A) circle (0.11);
                     \fill[black] (B) circle (0.11);
                     \node[above] at (A) {$\mathbf{p}_1$};
                     \node[below] at (B) {$\mathbf{p}_3$};
              \end{tikzpicture} +
              \begin{tikzpicture}[baseline = 1em, x = 1.4em, y = 1.4em]
                     \coordinate (A) at (0,0);
                     \coordinate (B) at (3,1.2);
                     \coordinate (C) at (5,2);
                     \draw[black] (B) -- (C);
                     \draw[dashed] (A) -- (B);
                     \fill[black] (B) circle (0.11);
                     \fill[black] (C) circle (0.11);
                     \node[below] at (C) {$\mathbf{p}_2$};
                     \node[below] at (B) {$\mathbf{p}_3$};
              \end{tikzpicture} -
              \begin{tikzpicture}[baseline = 1em, x = 1.4em, y = 1.4em]
                     \coordinate (A) at (0,0);
                     \coordinate (B) at (3,1.2);
                     \coordinate (C) at (5,2);
                     \draw[dashed] (A) -- (C);
                     \fill[black] (B) circle (0.11);
                     \node[below] at (B) {$\mathbf{p}_3$};
              \end{tikzpicture}.
       \end{gathered}
       \end{equation}
\end{lem}
\begin{proof}
       We have to prove that
       \begin{equation} \label{EquiOfPrfLine}
              \overline{\mathbf{p}_1 \mathbf{p}_2} + \mathbf{p}_3 = \overline{\mathbf{p}_1 \mathbf{p}_3} + \overline{\mathbf{p}_3 \mathbf{p}_2}.
       \end{equation}
       It holds that $\overline{\mathbf{p}_1 \mathbf{p}_2} + \mathbf{p}_3 = \overline{(\mathbf{p}_1 + \mathbf{p}_3) (\mathbf{p}_2 + \mathbf{p}_3)}$. On the other hand, $\overline{\mathbf{p}_1 \mathbf{p}_3} + \overline{\mathbf{p}_3 \mathbf{p}_2}$ can be rewritten as $\mathrm{conv}\{\mathbf{p}_1+\mathbf{p}_2,\mathbf{p}_1+\mathbf{p}_3,\mathbf{p}_2+\mathbf{p}_3,\mathbf{p}_3+\mathbf{p}_3\}$. The points $\mathbf{p}_1 + \mathbf{p}_2$ and $\mathbf{p}_3 + \mathbf{p}_3$ are inner points of $\overline{(\mathbf{p}_1 + \mathbf{p}_3) (\mathbf{p}_2 + \mathbf{p}_3)}$, since they can be expressed as 
       \begin{gather}
              \mathbf{p}_1 + \mathbf{p}_2= (1-t) (\mathbf{p}_1 + \mathbf{p}_3) + t (\mathbf{p}_2 + \mathbf{p}_3), \\
              \mathbf{p}_3 + \mathbf{p}_3 = t (\mathbf{p}_1 + \mathbf{p}_3) + (1-t) (\mathbf{p}_2 + \mathbf{p}_3),
       \end{gather}
       respectively. Therefore, $\mathrm{conv}\{\mathbf{p}_1+\mathbf{p}_2,\mathbf{p}_1+\mathbf{p}_3,\mathbf{p}_2+\mathbf{p}_3,\mathbf{p}_3+\mathbf{p}_3\} = \overline{(\mathbf{p}_1 + \mathbf{p}_3) (\mathbf{p}_2 + \mathbf{p}_3)}$ and \eqref{EquiOfPrfLine} holds.
\end{proof}
Lemma~\ref{LemOfLine} is applicable to segments. The following lemma is for polygons.
\begin{lem} \label{LemOfPolygon}
       Consider a polygon $\mathcal{P}$ which is separated into two polygons $\mathcal{P}_1$ and $\mathcal{P}_2$ by a straight line. Denote the intersection points between the boundary of $\mathcal{P}$ and the line as $\mathbf{l}_1$ and $\mathbf{l}_2$. Then it holds that
       \begin{equation} \label{EqOfLemOfPolygon}
       \begin{gathered}
              \mathcal{P} = \mathcal{P}_1 + \mathcal{P}_2 - \overline{\mathbf{l}_1 \mathbf{l}_2}. \\
              \Leftrightarrow \begin{tikzpicture}[baseline = 1.0em, x = 1.4em, y = 1.4em]
                     \coordinate (1) at (0.2,0);
                     \coordinate (2) at (1.4,-0.1);
                     \coordinate (3) at (2,0.8);
                     \coordinate (4) at (1.7,1.8);
                     \coordinate (5) at (0.4,2);
                     \coordinate (6) at (0,1);
                     \coordinate (7) at (0.8,-0.05);
                     \coordinate (8) at (0.7,1.94);
                     \fill[black!15] (1) -- (2) -- (3) -- (4) -- (5) -- (6) -- cycle;
                     \draw[black] (1) -- (2) -- (3) -- (4) -- (5) -- (6) -- cycle;
                     \draw[black,dashed] (7) -- (8);
                     \foreach \P in {1,2,3,4,5,6,7,8} \fill[black] (\P) circle (0.11);
              \end{tikzpicture}
              \,=\,
              \begin{tikzpicture}[baseline = 1.0em, x = 1.4em, y = 1.4em]
                     \coordinate (1) at (0.2,0);
                     \coordinate (2) at (1.4,-0.1);
                     \coordinate (3) at (2,0.8);
                     \coordinate (4) at (1.7,1.8);
                     \coordinate (5) at (0.4,2);
                     \coordinate (6) at (0,1);
                     \coordinate (7) at (0.8,-0.05);
                     \coordinate (8) at (0.7,1.94);
                     \fill[black!15] (1) -- (7) -- (8) -- (5) -- (6) -- cycle;
                     \draw[black] (1) -- (7) -- (8) -- (5) -- (6) -- cycle;
                     \draw[black, dashed] (7) -- (2) -- (3) -- (4) -- (8) -- cycle;
                     \foreach \P in {1,5,6,7,8} \fill[black] (\P) circle (0.11);
              \end{tikzpicture} \,+\,
              \begin{tikzpicture}[baseline = 1.0em, x = 1.4em, y = 1.4em]
                     \coordinate (1) at (0.2,0);
                     \coordinate (2) at (1.4,-0.1);
                     \coordinate (3) at (2,0.8);
                     \coordinate (4) at (1.7,1.8);
                     \coordinate (5) at (0.4,2);
                     \coordinate (6) at (0,1);
                     \coordinate (7) at (0.8,-0.05);
                     \coordinate (8) at (0.7,1.94);
                     \fill[black!15] (7) -- (2) -- (3) -- (4) -- (8) -- cycle;
                     \draw[black, dashed] (1) -- (7) -- (8) -- (5) -- (6) -- cycle;
                     \draw[black] (7) -- (2) -- (3) -- (4) -- (8) -- cycle;
                     \foreach \P in {2,3,4,7,8} \fill[black] (\P) circle (0.11);
              \end{tikzpicture} \,-\,
              \begin{tikzpicture}[baseline = 1.0em, x = 1.4em, y = 1.4em]
                     \coordinate (1) at (0.2,0);
                     \coordinate (2) at (1.4,-0.1);
                     \coordinate (3) at (2,0.8);
                     \coordinate (4) at (1.7,1.8);
                     \coordinate (5) at (0.4,2);
                     \coordinate (6) at (0,1);
                     \coordinate (7) at (0.8,-0.05);
                     \coordinate (8) at (0.7,1.94);
                     \draw[black, dashed] (1) -- (2) -- (3) -- (4) -- (5) -- (6) -- cycle;
                     \draw[black] (7) -- (8);
                     \foreach \P in {7,8} \fill[black] (\P) circle (0.11);
              \end{tikzpicture}\,.
       \end{gathered}
\end{equation}
\end{lem}
\begin{proof}
       For proving \eqref{EqOfLemOfPolygon}, it is sufficient to verify
\begin{equation} \label{EqOfProveOfLemOfPolygon}
       \mathcal{P} + \overline{\mathbf{l}_1 \mathbf{l}_2} = \mathcal{P}_1 + \mathcal{P}_2.
\end{equation}
The left-hand side of \eqref{EqOfProveOfLemOfPolygon} can be rewritten as
\begin{equation} \label{EqOfProveOfLemOfPolygonleft}
       \begin{aligned}
              \mathcal{P} + \overline{\mathbf{l}_1 \mathbf{l}_2} &= \mathrm{conv}(v(\mathcal{P})) + \overline{\mathbf{l}_1 \mathbf{l}_2}\\
              &= \mathrm{conv}(v(\mathcal{P}) \cup \{\mathbf{l}_1, \mathbf{l}_2\}) + \overline{\mathbf{l}_1 \mathbf{l}_2} \\
              &= \mathrm{conv}((v(\mathcal{P}) \cup \{\mathbf{l}_1, \mathbf{l}_2\}) + \{\mathbf{l}_1, \mathbf{l}_2\}) \\
              &= \mathrm{conv}\{\mathbf{l}_1 + \mathbf{l}_2, 2\mathbf{l}_1, 2\mathbf{l}_2, \mathbf{p}_i + \mathbf{l}_1, \mathbf{p}_i + \mathbf{l}_2 \,|\, \mathbf{p}_i \in v(\mathcal{P})\} \\
              &= \mathrm{conv}\{\mathbf{p}_i + \mathbf{l}_1, \mathbf{p}_i + \mathbf{l}_2 \,|\, \mathbf{p}_i \in v(\mathcal{P})\}
       \end{aligned}
\end{equation}
according to Lemma~\ref{LemOfGeometricIdentities}. Note that the terms $\mathbf{l}_1 + \mathbf{l}_2$, $2\mathbf{l}_1$, and $2\mathbf{l}_2$ can be removed, since $\mathbf{l}_1$ and $\mathbf{l}_2$ are vertices of $\mathcal{P}$ or interior points of an edge of $\mathcal{P}$. The right-hand side of \eqref{EqOfProveOfLemOfPolygon} can be rewritten as
\begin{equation} \label{EqOfProveOfLemOfPolygonright}
       \begin{aligned}
              \mathcal{P}_1 + \mathcal{P}_2 = &\mathrm{conv}(v(\mathcal{P}_1)) + \mathrm{conv}(v(\mathcal{P}_2)) \\
              = &\mathrm{conv}\{\mathbf{p}_i + \mathbf{l}_1, \mathbf{p}_i + \mathbf{l}_2, \mathbf{p}_j + \mathbf{p}_k \,|\, \mathbf{p}_i \in v(\mathcal{P}),\, \mathbf{p}_j \in v(\mathcal{P}_1),\, \mathbf{p}_k \in v(\mathcal{P}_2)\}.
       \end{aligned}
\end{equation}
The difference between \eqref{EqOfProveOfLemOfPolygonleft} and \eqref{EqOfProveOfLemOfPolygonright} is given by the terms $\mathbf{p}_j + \mathbf{p}_k$. We will show that, for any pair of $\mathbf{p}_j \in v(\mathcal{P}_1)$ and $\mathbf{p}_k \in v(\mathcal{P}_2)$, the point $\mathbf{p}_j + \mathbf{p}_k$ is contained in $\mathrm{conv}\{\mathbf{p}_j + \mathbf{l}_1, \mathbf{p}_k + \mathbf{l}_1, \mathbf{p}_j + \mathbf{l}_2, \mathbf{p}_k + \mathbf{l}_2\}$ as shown in Figure~\ref{FigOfProofPolygonCut}. If $\mathbf{p}_j$ or $\mathbf{p}_k$ correspond to $\mathbf{l}_1$ or $\mathbf{l}_2$, $\mathbf{p}_j + \mathbf{p}_k$ is one of $\mathbf{p}_j + \mathbf{l}_1$, $\mathbf{p}_j + \mathbf{l}_2$, $\mathbf{p}_k + \mathbf{l}_1$, or $\mathbf{p}_k + \mathbf{l}_2$. Therefore, $\mathbf{p}_j + \mathbf{p}_k$ is contained in $\mathrm{conv}\{\mathbf{p}_j + \mathbf{l}_1, \mathbf{p}_k + \mathbf{l}_1, \mathbf{p}_j + \mathbf{l}_2, \mathbf{p}_k + \mathbf{l}_2\}$. Consider the cases such that $\mathbf{p}_j$ and $\mathbf{p}_k$ are different from $\mathbf{l}_1$ and $\mathbf{l}_2$. Fix $\mathbf{p}_j = (a_j,b_j), \mathbf{p}_k = (a_k,b_k), \mathbf{l}_1 = (a_1,b_1), \mathbf{l}_2 = (a_2,b_2)$, and define $f_{jk} (x,y) := (b_j - b_k)x - (a_j - a_k)y + a_j b_k - a_k b_j$ and $f_{l} (x,y) := (b_1 - b_2)x - (a_1 - a_2)y + a_1 b_2 - a_2 b_1$ respectively. Note that $f_{jk}(\mathbf{p}_j) = f_{jk}(\mathbf{p}_k) = f_{l}(\mathbf{l}_1) = f_{l}(\mathbf{l}_2) = 0$, and the relative positions between a point $\mathbf{p}$ and $\overline{\mathbf{p}_j\mathbf{p}_k}$ is determined by whether $f_{jk}(\mathbf{p})$ is positive, $0$, or negative. Similarly, the relative positions between $\mathbf{p}$ and $\overline{\mathbf{l}_1 \mathbf{l}_2}$ is determined by the value of $f_{l}(\mathbf{p})$. Without loss of generality, we can assume
\begin{equation}
       \begin{cases}
              f_{jk}(\mathbf{l}_1) > 0, \\
              f_{jk}(\mathbf{l}_2) < 0, \\
              f_{l}(\mathbf{p}_j) > 0, \\
              f_{l}(\mathbf{p}_k) < 0.
       \end{cases}
\end{equation}
From these inequalities and linearity of $f_{jk}$ and $f_{l}$, we can obtain
\begin{equation}
\begin{cases}
       f_{jk}(\mathbf{p}_j + \mathbf{p}_k - \mathbf{l}_1) = - f_{jk}(\mathbf{l}_1) < 0, \\
       f_{jk}(\mathbf{p}_j + \mathbf{p}_k - \mathbf{l}_2) = - f_{jk}(\mathbf{l}_2) > 0, \\
       f_{l}(\mathbf{p}_j + \mathbf{p}_k - \mathbf{p}_j) = f_{l}(\mathbf{p}_k) < 0, \\
       f_{l}(\mathbf{p}_j + \mathbf{p}_k - \mathbf{p}_k) = f_{l}(\mathbf{p}_j) > 0.
\end{cases}
\end{equation}
These inequalities indicate that the positions of $\mathbf{p}_j + \mathbf{p}_k, \mathbf{p}_j + \mathbf{l}_1, \mathbf{p}_k + \mathbf{l}_1, \mathbf{p}_j + \mathbf{l}_2, \mathbf{p}_k + \mathbf{l}_2$ are determined as shown in Figure~\ref{FigOfProofPolygonCut}. 
Therefore $\mathbf{p}_j + \mathbf{p}_k$ must belong to $\mathrm{conv}\{\mathbf{p}_j + \mathbf{l}_1, \mathbf{p}_j + \mathbf{l}_2, \mathbf{p}_k + \mathbf{l}_1, \mathbf{p}_k + \mathbf{l}_2\}$ for any pair of $\mathbf{p}_j$ and $\mathbf{p}_k$, and we can remove the terms $\mathbf{p}_j + \mathbf{p}_k$ from \eqref{EqOfProveOfLemOfPolygonright}. As a result, \eqref{EqOfProveOfLemOfPolygon} follows.
\end{proof}
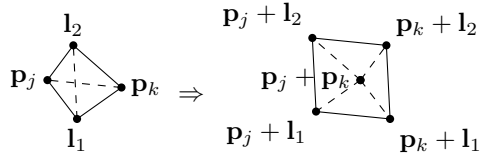
\begin{figure}
       \centering
       \begin{tikzpicture}[baseline = 0.6em, x = 1.4em, y = 1.4em]
              \coordinate (1) at (0.2, 0);
              \coordinate (2) at (1.4, -0.1);
              \coordinate (3) at (2, 0.8);
              \coordinate (4) at (1.7, 1.8);
              \coordinate (5) at (0.4, 2);
              \coordinate (6) at (0, 1);
              \coordinate (7) at (0.8, -0.05);
              \coordinate (8) at (0.7, 1.94);
              \draw[black] (3) -- (8) -- (6) -- (7) -- cycle;
              \draw[dashed] (7) -- (8);
              \draw[dashed] (3) -- (6);
              \foreach \P in {3,6,7,8} \fill[black] (\P) circle (0.11);
              \node[right] at (3) {$\mathbf{p}_k$};
              \node[left] at (6) {$\mathbf{p}_j$};
              \node[below] at (7) {$\mathbf{l}_1$};
              \node[above] at (8) {$\mathbf{l}_2$};
       \end{tikzpicture} \text{$\Rightarrow$}
       \begin{tikzpicture}[baseline = 0.6em, x = 1.4em, y = 1.4em]
              \coordinate (0) at (0.4, 2);
              \coordinate (1) at (0, 1);
              \coordinate (2) at (0.2, 0);
              \coordinate (3) at (0.8, -0.05);
              \coordinate (4) at (0.7, 1.94);
              \coordinate (5) at (1.4 - 2, -0.1 + 0.2);
              \coordinate (6) at (2 - 2, 0.8 + 0.2);
              \coordinate (7) at (1.7 - 2, 1.8 + 0.2);
              \coordinate (8) at (0.7 - 2, 1.94 + 0.2);
              \coordinate (9) at (0.8 - 2, -0.05 + 0.2);
              \draw[black] (3) -- (4);
              \draw[black] (9) -- (8);
              \draw[black] (3) -- (9);
              \draw[black] (4) -- (8);
              \draw[black, dashed] (1) -- (3);
              \draw[black, dashed] (1) -- (4);
              \draw[black, dashed] (1) -- (8);
              \draw[black, dashed] (1) -- (9);
              \node[below right] at (3) {$\mathbf{p}_k + \mathbf{l}_1$};
              \node[above right] at (4) {$\mathbf{p}_k + \mathbf{l}_2$};
              \node[above left] at (8) {$\mathbf{p}_j + \mathbf{l}_2$};
              \node[below left] at (9) {$\mathbf{p}_j + \mathbf{l}_1$};
              \node[left] at (1) {$\mathbf{p}_j + \mathbf{p}_k$};
              \foreach \P in {1,3,4,6,8,9} \fill[black] (\P) circle (0.11);
       \end{tikzpicture}
       \caption{Figure for the proof of Lemma~\ref{LemOfPolygon}}
       \label{FigOfProofPolygonCut}
\end{figure}
\subsection{Signed Minkowski decompositions induced from partitions of polygons}
For a segment or a polygon $\mathcal{C}$, if it has an expression as a Minkowski sum of points, segments, or polygons, $\mathcal{C}_1$, $\dots$, $\mathcal{C}_m$, namely
\begin{equation}
       \mathcal{C} = \sum_{i=1}^{m} \mathcal{C}_i := \mathcal{C}_1 + \dots + \mathcal{C}_m, \label{EqOfMinkowskiDecomp}
\end{equation}
where at least two of $\mathcal{C}_i$ are either a segment or a polygon, then the right-hand side of \eqref{EqOfMinkowskiDecomp} is called a Minkowski decomposition of $\mathcal{C}$. If $\mathcal{C}$ is expressed by Minkowski sums and differences,
\begin{equation}
       \begin{aligned}
              \mathcal{C} &= \sum_{i=1}^{m} \mathcal{C}_i' - \sum_{j=1}^{n} \mathcal{C}_j'' \\
              &= \mathcal{C}_1' + \dots + \mathcal{C}_m' - \mathcal{C}_1'' - \dots - \mathcal{C}_n'', \label{EqOfSignedMinkowskiDecomp}
       \end{aligned}
\end{equation}
then the right-hand side of \eqref{EqOfSignedMinkowskiDecomp} is called a signed Minkowski decomposition of $\mathcal{C}$.

A partition of a polygon is its subdivision into smaller polygons~\cite{Khovanskii1997}. Funke elucidated that if a partition of a polygon is given, a signed Minkowski decomposition of the polygon can be induced by applying the cutting relations along dividing lines of the partition~\cite{Funke2021}. This decomposition is referred as a `Partition relation'~\cite[Proposition 3.8]{Funke2021}. Before deriving decompositions induced by general partitions, we explain how to derive a decomposition via a simple example. For convenience, we employ the notation $\mathrm{conv} \{\mathbf{p}_1, \mathbf{p}_2, \mathbf{p}_3\} = \triangle \mathbf{p}_1 \mathbf{p}_2 \mathbf{p}_3$ if $\mathrm{conv}\{\mathbf{p}_1$, $\mathbf{p}_2$, $\mathbf{p}_3\}$ is a triangle. 
\begin{cor} \label{CorTriangleSubdivision}
       Assume that $\mathbf{p}_4$ is an interior point of $\triangle \mathbf{p}_1 \mathbf{p}_2 \mathbf{p}_3$. Then $\triangle \mathbf{p}_1 \mathbf{p}_2 \mathbf{p}_3$ has a signed decomposition expressed as
       \begin{equation} \label{EqCorTriangleSubdivision}
              \begin{aligned}
                     \triangle \mathbf{p}_1 \mathbf{p}_2 \mathbf{p}_3 = &\triangle \mathbf{p}_1 \mathbf{p}_2 \mathbf{p}_4 + \triangle \mathbf{p}_2 \mathbf{p}_3 \mathbf{p}_4 + \triangle \mathbf{p}_1 \mathbf{p}_3 \mathbf{p}_4 \\
                     &- \overline{\mathbf{p}_1 \mathbf{p}_4} - \overline{\mathbf{p}_2 \mathbf{p}_4} - \overline{\mathbf{p}_3 \mathbf{p}_4} + \mathbf{p}_4.
              \end{aligned}
       \end{equation}
\end{cor}
\begin{rem}
The decomposition \eqref{EqCorTriangleSubdivision} is visualized as follows,
\begin{equation}
\begin{aligned}
       \begin{tikzpicture}[baseline = 1.0em, x = 1.4em, y = 1.4em]
              \coordinate (1) at (0, 0);
              \coordinate (2) at (3, 0);
              \coordinate (3) at (2, 2);
              \coordinate (4) at (1.5, 0.7);
              \fill[black!15] (1) -- (2) -- (3) -- cycle;
              \draw (1) -- (2) -- (3) -- cycle;
              \draw[dashed] (1) -- (4);
              \draw[dashed] (2) -- (4);
              \draw[dashed] (3) -- (4);
              \foreach \P in {1,2,3,4} \fill[black] (\P) circle (0.11);
              \node[left] at (1) {$\mathbf{p}_1$};
              \node[right] at (2) {$\mathbf{p}_2$};
              \node[above] at (3) {$\mathbf{p}_3$};
              \node[below] at (4) {$\mathbf{p}_4$};
       \end{tikzpicture} =& 
       \begin{tikzpicture}[baseline = 1.0em, x = 1.4em, y = 1.4em]
              \coordinate (1) at (0, 0);
              \coordinate (2) at (3, 0);
              \coordinate (3) at (2, 2);
              \coordinate (4) at (1.5, 0.7);
              \fill[black!15] (1) -- (2) -- (4) -- cycle;
              \draw[dashed] (1) -- (2) -- (3) -- cycle;
              \draw (1) -- (2) -- (4) -- cycle;
              \draw[dashed] (1) -- (4);
              \draw[dashed] (2) -- (4);
              \draw[dashed] (3) -- (4);
              \foreach \P in {1,2,4} \fill[black] (\P) circle (0.11);
              \node[left] at (1) {$\mathbf{p}_1$};
              \node[right] at (2) {$\mathbf{p}_2$};
              \node[below] at (4) {$\mathbf{p}_4$};
       \end{tikzpicture} +
       \begin{tikzpicture}[baseline = 1.0em, x = 1.4em, y = 1.4em]
              \coordinate (1) at (0, 0);
              \coordinate (2) at (3, 0);
              \coordinate (3) at (2, 2);
              \coordinate (4) at (1.5, 0.7);
              \fill[black!15] (4) -- (2) -- (3) -- cycle;
              \draw (4) -- (2) -- (3) -- cycle;
              \draw[dashed] (1) -- (2) -- (3) -- cycle;
              \draw[dashed] (1) -- (4);
              \draw[dashed] (2) -- (4);
              \draw[dashed] (3) -- (4);
              \foreach \P in {2,3,4} \fill[black] (\P) circle (0.11);
              \node[right] at (2) {$\mathbf{p}_2$};
              \node[above] at (3) {$\mathbf{p}_3$};
              \node[below] at (4) {$\mathbf{p}_4$};
       \end{tikzpicture} +
       \begin{tikzpicture}[baseline = 1.0em, x = 1.4em, y = 1.4em]
              \coordinate (1) at (0, 0);
              \coordinate (2) at (3, 0);
              \coordinate (3) at (2, 2);
              \coordinate (4) at (1.5, 0.7);
              \fill[black!15] (1) -- (4) -- (3) -- cycle;
              \draw (1) -- (4) -- (3) -- cycle;
              \draw[dashed] (1) -- (2) -- (3) -- cycle;
              \draw[dashed] (1) -- (4);
              \draw[dashed] (2) -- (4);
              \draw[dashed] (3) -- (4);
              \foreach \P in {1,3,4} \fill[black] (\P) circle (0.11);
              \node[left] at (1) {$\mathbf{p}_1$};
              \node[above] at (3) {$\mathbf{p}_3$};
              \node[below] at (4) {$\mathbf{p}_4$};
       \end{tikzpicture} \\
       &-
       \begin{tikzpicture}[baseline = 1.0em, x = 1.4em, y = 1.4em]
              \coordinate (1) at (0, 0);
              \coordinate (2) at (3, 0);
              \coordinate (3) at (2, 2);
              \coordinate (4) at (1.5, 0.7);
              \draw[dashed] (1) -- (2) -- (3) -- cycle;
              \draw[thick] (1) -- (4);
              \draw[dashed] (2) -- (4);
              \draw[dashed] (3) -- (4);
              \foreach \P in {1,4} \fill[black] (\P) circle (0.11);
              \node[left] at (1) {$\mathbf{p}_1$};
              \node[below] at (4) {$\mathbf{p}_4$};
       \end{tikzpicture} -
       \begin{tikzpicture}[baseline = 1.0em, x = 1.4em, y = 1.4em]
              \coordinate (1) at (0, 0);
              \coordinate (2) at (3, 0);
              \coordinate (3) at (2, 2);
              \coordinate (4) at (1.5, 0.7);
              \draw[dashed] (1) -- (2) -- (3) -- cycle;
              \draw[dashed] (1) -- (4);
              \draw[thick] (2) -- (4);
              \draw[dashed] (3) -- (4);
              \foreach \P in {2,4} \fill[black] (\P) circle (0.11);
              \node[right] at (2) {$\mathbf{p}_2$};
              \node[below] at (4) {$\mathbf{p}_4$};
       \end{tikzpicture} -
       \begin{tikzpicture}[baseline = 1.0em, x = 1.4em, y = 1.4em]
              \coordinate (1) at (0, 0);
              \coordinate (2) at (3, 0);
              \coordinate (3) at (2, 2);
              \coordinate (4) at (1.5, 0.7);
              \draw[dashed] (1) -- (2) -- (3) -- cycle;
              \draw[dashed] (1) -- (4);
              \draw[dashed] (2) -- (4);
              \draw[thick] (3) -- (4);
              \foreach \P in {3,4} \fill[black] (\P) circle (0.11);
              \node[above] at (3) {$\mathbf{p}_3$};
              \node[below] at (4) {$\mathbf{p}_4$};
       \end{tikzpicture} +
       \begin{tikzpicture}[baseline = 1.0em, x = 1.4em, y = 1.4em]
              \coordinate (1) at (0, 0);
              \coordinate (2) at (3, 0);
              \coordinate (3) at (2, 2);
              \coordinate (4) at (1.5, 0.7);
              \draw[dashed] (1) -- (2) -- (3) -- cycle;
              \draw[dashed] (1) -- (4);
              \draw[dashed] (2) -- (4);
              \draw[dashed] (3) -- (4);
              \foreach \P in {4} \fill[black] (\P) circle (0.11);
              \node[below] at (4) {$\mathbf{p}_4$};
       \end{tikzpicture}.
\end{aligned}
\end{equation}
\end{rem}
\begin{proof}
       Firstly, we define $\mathbf{p}_5$ as the intersection point of the extension of $\overline{\mathbf{p}_1 \mathbf{p}_4}$ and $\overline{\mathbf{p}_2 \mathbf{p}_3}$. Applying Lemma~\ref{LemOfPolygon} to $\overline{\mathbf{p}_3 \mathbf{p}_5}$, we obtain
       \begin{equation}
              \triangle \mathbf{p}_1 \mathbf{p}_2 \mathbf{p}_3 = \triangle \mathbf{p}_1 \mathbf{p}_3 \mathbf{p}_5 + \triangle \mathbf{p}_2 \mathbf{p}_3 \mathbf{p}_5 - \overline{\mathbf{p}_3 \mathbf{p}_5}. 
       \end{equation}
       Similarly, by applying Lemma~\ref{LemOfLine} to $\overline{\mathbf{p}_3 \mathbf{p}_5}$ and Lemma~\ref{LemOfPolygon} to $\triangle \mathbf{p}_1 \mathbf{p}_3 \mathbf{p}_5$ and $\triangle \mathbf{p}_2 \mathbf{p}_3 \mathbf{p}_5$, we obtain
       \begin{equation} \label{Eq1PrfCorTriangleSubdivision}
       \begin{aligned}
              \triangle \mathbf{p}_1 \mathbf{p}_3 \mathbf{p}_5 + \triangle \mathbf{p}_2 \mathbf{p}_3 \mathbf{p}_5 - \overline{\mathbf{p}_3 \mathbf{p}_5}
              =& \triangle \mathbf{p}_1 \mathbf{p}_3 \mathbf{p}_4 + \triangle \mathbf{p}_1 \mathbf{p}_4 \mathbf{p}_5 - \overline{\mathbf{p}_1 \mathbf{p}_4} \\
              &+ \triangle \mathbf{p}_2 \mathbf{p}_3 \mathbf{p}_4 + \triangle \mathbf{p}_2 \mathbf{p}_4 \mathbf{p}_5 - \overline{\mathbf{p}_2 \mathbf{p}_4} \\
              &- \overline{\mathbf{p}_3 \mathbf{p}_4} - \overline{\mathbf{p}_4 \mathbf{p}_5} + \mathbf{p}_4.
       \end{aligned}
       \end{equation}
       By applying Lemma~\ref{LemOfPolygon}, it holds that
       \begin{equation} \label{Eq2PrfCorTriangleSubdivision}
              \triangle \mathbf{p}_1 \mathbf{p}_4 \mathbf{p}_5 + \triangle \mathbf{p}_2 \mathbf{p}_4 \mathbf{p}_5 - \overline{\mathbf{p}_4 \mathbf{p}_5} = \triangle \mathbf{p}_1 \mathbf{p}_2 \mathbf{p}_4.
       \end{equation}
       By substituting \eqref{Eq2PrfCorTriangleSubdivision} into \eqref{Eq1PrfCorTriangleSubdivision}, we obtain \eqref{EqCorTriangleSubdivision}.
\end{proof}
\begin{rem}
       The visualized proof of Corollary~\ref{CorTriangleSubdivision} is as follows,
       \begin{equation}
       \begin{aligned}
              \begin{tikzpicture}[baseline = 1.3em, x = 1.4em, y = 1.4em]
                     \coordinate (1) at (0, 0);
                     \coordinate (2) at (3, 0);
                     \coordinate (3) at (2, 2);
                     \coordinate (4) at (1.5, 0.7);
                     \coordinate (5) at (2 - 10 / 13, 0);
                     \fill[black!15] (1) -- (2) -- (3) -- cycle;
                     \draw (1) -- (2) -- (3) -- cycle;
                     \draw[dashed] (3) -- (5);
                     \foreach \P in {1,2,3,4,5} \fill[black] (\P) circle (0.11);
                     \node[left] at (1) {$\mathbf{p}_1$};
                     \node[right] at (2) {$\mathbf{p}_2$};
                     \node[above] at (3) {$\mathbf{p}_3$};
                     \node[left] at (4) {$\mathbf{p}_4$};
                     \node[below] at (5) {$\mathbf{p}_5$};
              \end{tikzpicture} = &
              \begin{tikzpicture}[baseline = 1.3em, x = 1.4em, y = 1.4em]
                     \coordinate (1) at (0, 0);
                     \coordinate (2) at (3, 0);
                     \coordinate (3) at (2, 2);
                     \coordinate (4) at (1.5, 0.7);
                     \coordinate (5) at (2 - 10 / 13, 0);
                     \fill[black!15] (1) -- (5) -- (3) -- cycle;
                     \draw (1) -- (5) -- (3) -- cycle;
                     \draw[dashed] (1) -- (2) -- (3) -- cycle;
                     \draw[dashed] (1) -- (4);
                     \foreach \P in {1,3,4,5} \fill[black] (\P) circle (0.11);
                     \node[left] at (1) {$\mathbf{p}_1$};
                     \node[above] at (3) {$\mathbf{p}_3$};
                     \node[left] at (4) {$\mathbf{p}_4$};
                     \node[below] at (5) {$\mathbf{p}_5$};
              \end{tikzpicture} +
              \begin{tikzpicture}[baseline = 1.3em, x = 1.4em, y = 1.4em]
                     \coordinate (1) at (0, 0);
                     \coordinate (2) at (3, 0);
                     \coordinate (3) at (2, 2);
                     \coordinate (4) at (1.5, 0.7);
                     \coordinate (5) at (2 - 10 / 13, 0);
                     \fill[black!15] (2) -- (5) -- (3) -- cycle;
                     \draw (2) -- (5) -- (3) -- cycle;
                     \draw[dashed] (1) -- (2) -- (3) -- cycle;
                     \draw[dashed] (2) -- (4);
                     \foreach \P in {2,3,4,5} \fill[black] (\P) circle (0.11);
                     \node[right] at (2) {$\mathbf{p}_2$};
                     \node[above] at (3) {$\mathbf{p}_3$};
                     \node[left] at (4) {$\mathbf{p}_4$};
                     \node[below] at (5) {$\mathbf{p}_5$};
              \end{tikzpicture} -
              \begin{tikzpicture}[baseline = 1.3em, x = 1.4em, y = 1.4em]
                     \coordinate (1) at (0, 0);
                     \coordinate (2) at (3, 0);
                     \coordinate (3) at (2, 2);
                     \coordinate (4) at (1.5, 0.7);
                     \coordinate (5) at (2 - 10 / 13, 0);
                     \draw[dashed] (1) -- (2) -- (3) -- cycle;
                     \draw (3) -- (5);
                     \foreach \P in {3,4,5} \fill[black] (\P) circle (0.11);
                     \node[above] at (3) {$\mathbf{p}_3$};
                     \node[left] at (4) {$\mathbf{p}_4$};
                     \node[below] at (5) {$\mathbf{p}_5$};
              \end{tikzpicture} \\
              = &
              \begin{tikzpicture}[baseline = 1.3em, x = 1.4em, y = 1.4em]
                     \coordinate (1) at (0, 0);
                     \coordinate (2) at (3, 0);
                     \coordinate (3) at (2, 2);
                     \coordinate (4) at (1.5, 0.7);
                     \coordinate (5) at (2 - 10 / 13, 0);
                     \fill[black!15] (1) -- (4) -- (3) -- cycle;
                     \draw (1) -- (4) -- (3) -- cycle;
                     \draw[dashed] (1) -- (2) -- (3) -- cycle;
                     \foreach \P in {1,3,4} \fill[black] (\P) circle (0.11);
                     \node[left] at (1) {$\mathbf{p}_1$};
                     \node[above] at (3) {$\mathbf{p}_3$};
                     \node[left] at (4) {$\mathbf{p}_4$};
              \end{tikzpicture} +
              \begin{tikzpicture}[baseline = 1.3em, x = 1.4em, y = 1.4em]
                     \coordinate (1) at (0, 0);
                     \coordinate (2) at (3, 0);
                     \coordinate (3) at (2, 2);
                     \coordinate (4) at (1.5, 0.7);
                     \coordinate (5) at (2 - 10 / 13, 0);
                     \fill[black!15] (1) -- (5) -- (4) -- cycle;
                     \draw (1) -- (5) -- (4) -- cycle;
                     \draw[dashed] (1) -- (2) -- (3) -- cycle;
                     \foreach \P in {1,4,5} \fill[black] (\P) circle (0.11);
                     \node[left] at (1) {$\mathbf{p}_1$};
                     \node[left] at (4) {$\mathbf{p}_4$};
                     \node[below] at (5) {$\mathbf{p}_5$};
              \end{tikzpicture} -
              \begin{tikzpicture}[baseline = 1.3em, x = 1.4em, y = 1.4em]
                     \coordinate (1) at (0, 0);
                     \coordinate (2) at (3, 0);
                     \coordinate (3) at (2, 2);
                     \coordinate (4) at (1.5, 0.7);
                     \coordinate (5) at (2 - 10 / 13, 0);
                     \draw[dashed] (1) -- (2) -- (3) -- cycle;
                     \draw (1) -- (4);
                     \foreach \P in {1,4} \fill[black] (\P) circle (0.11);
                     \node[left] at (1) {$\mathbf{p}_1$};
                     \node[left] at (4) {$\mathbf{p}_4$};
              \end{tikzpicture} \\
              &+
              \begin{tikzpicture}[baseline = 1.3em, x = 1.4em, y = 1.4em]
                     \coordinate (1) at (0, 0);
                     \coordinate (2) at (3, 0);
                     \coordinate (3) at (2, 2);
                     \coordinate (4) at (1.5, 0.7);
                     \coordinate (5) at (2 - 10 / 13, 0);
                     \fill[black!15] (2) -- (4) -- (3) -- cycle;
                     \draw (2) -- (4) -- (3) -- cycle;
                     \draw[dashed] (1) -- (2) -- (3) -- cycle;
                     \foreach \P in {2,3,4} \fill[black] (\P) circle (0.11);
                     \node[right] at (2) {$\mathbf{p}_2$};
                     \node[above] at (3) {$\mathbf{p}_3$};
                     \node[left] at (4) {$\mathbf{p}_4$};
              \end{tikzpicture} +
              \begin{tikzpicture}[baseline = 1.3em, x = 1.4em, y = 1.4em]
                     \coordinate (1) at (0, 0);
                     \coordinate (2) at (3, 0);
                     \coordinate (3) at (2, 2);
                     \coordinate (4) at (1.5, 0.7);
                     \coordinate (5) at (2 - 10 / 13, 0);
                     \fill[black!15] (2) -- (5) -- (4) -- cycle;
                     \draw (2) -- (5) -- (4) -- cycle;
                     \draw[dashed] (1) -- (2) -- (3) -- cycle;
                     \foreach \P in {2,4,5} \fill[black] (\P) circle (0.11);
                     \node[right] at (2) {$\mathbf{p}_2$};
                     \node[left] at (4) {$\mathbf{p}_4$};
                     \node[below] at (5) {$\mathbf{p}_5$};
              \end{tikzpicture} -
              \begin{tikzpicture}[baseline = 1.3em, x = 1.4em, y = 1.4em]
                     \coordinate (1) at (0, 0);
                     \coordinate (2) at (3, 0);
                     \coordinate (3) at (2, 2);
                     \coordinate (4) at (1.5, 0.7);
                     \coordinate (5) at (2 - 10 / 13, 0);
                     \draw (2) -- (4) -- cycle;
                     \draw[dashed] (1) -- (2) -- (3) -- cycle;
                     \foreach \P in {2,4} \fill[black] (\P) circle (0.11);
                     \node[right] at (2) {$\mathbf{p}_2$};
                     \node[left] at (4) {$\mathbf{p}_4$};
              \end{tikzpicture} \\
              &-
              \begin{tikzpicture}[baseline = 1.3em, x = 1.4em, y = 1.4em]
                     \coordinate (1) at (0, 0);
                     \coordinate (2) at (3, 0);
                     \coordinate (3) at (2, 2);
                     \coordinate (4) at (1.5, 0.7);
                     \coordinate (5) at (2 - 10 / 13, 0);
                     \draw[dashed] (1) -- (2) -- (3) -- cycle;
                     \draw (3) -- (4);
                     \foreach \P in {3,4} \fill[black] (\P) circle (0.11);
                     \node[above] at (3) {$\mathbf{p}_3$};
                     \node[left] at (4) {$\mathbf{p}_4$};
              \end{tikzpicture} -
              \begin{tikzpicture}[baseline = 1.3em, x = 1.4em, y = 1.4em]
                     \coordinate (1) at (0, 0);
                     \coordinate (2) at (3, 0);
                     \coordinate (3) at (2, 2);
                     \coordinate (4) at (1.5, 0.7);
                     \coordinate (5) at (2 - 10 / 13, 0);
                     \draw[dashed] (1) -- (2) -- (3) -- cycle;
                     \draw (4) -- (5);
                     \foreach \P in {4,5} \fill[black] (\P) circle (0.11);
                     \node[left] at (4) {$\mathbf{p}_4$};
                     \node[below] at (5) {$\mathbf{p}_5$};
              \end{tikzpicture} +
              \begin{tikzpicture}[baseline = 1.3em, x = 1.4em, y = 1.4em]
                     \coordinate (1) at (0, 0);
                     \coordinate (2) at (3, 0);
                     \coordinate (3) at (2, 2);
                     \coordinate (4) at (1.5, 0.7);
                     \coordinate (5) at (2 - 10 / 13, 0);
                     \draw[dashed] (1) -- (2) -- (3) -- cycle;
                     \foreach \P in {4} \fill[black] (\P) circle (0.11);
                     \node[left] at (4) {$\mathbf{p}_4$};
              \end{tikzpicture} \\
              =&
              \begin{tikzpicture}[baseline = 1.0em, x = 1.4em, y = 1.4em]
                     \coordinate (1) at (0, 0);
                     \coordinate (2) at (3, 0);
                     \coordinate (3) at (2, 2);
                     \coordinate (4) at (1.5, 0.7);
                     \fill[black!15] (1) -- (2) -- (4) -- cycle;
                     \draw[dashed] (1) -- (2) -- (3) -- cycle;
                     \draw (1) -- (2) -- (4) -- cycle;
                     \draw[dashed] (1) -- (4);
                     \draw[dashed] (2) -- (4);
                     \draw[dashed] (3) -- (4);
                     \foreach \P in {1,2,4} \fill[black] (\P) circle (0.11);
                     \node[left] at (1) {$\mathbf{p}_1$};
                     \node[right] at (2) {$\mathbf{p}_2$};
                     \node[below] at (4) {$\mathbf{p}_4$};
              \end{tikzpicture} +
              \begin{tikzpicture}[baseline = 1.0em, x = 1.4em, y = 1.4em]
                     \coordinate (1) at (0, 0);
                     \coordinate (2) at (3, 0);
                     \coordinate (3) at (2, 2);
                     \coordinate (4) at (1.5, 0.7);
                     \fill[black!15] (4) -- (2) -- (3) -- cycle;
                     \draw (4) -- (2) -- (3) -- cycle;
                     \draw[dashed] (1) -- (2) -- (3) -- cycle;
                     \draw[dashed] (1) -- (4);
                     \draw[dashed] (2) -- (4);
                     \draw[dashed] (3) -- (4);
                     \foreach \P in {2,3,4} \fill[black] (\P) circle (0.11);
                     \node[right] at (2) {$\mathbf{p}_2$};
                     \node[above] at (3) {$\mathbf{p}_3$};
                     \node[below] at (4) {$\mathbf{p}_4$};
              \end{tikzpicture} +
              \begin{tikzpicture}[baseline = 1.0em, x = 1.4em, y = 1.4em]
                     \coordinate (1) at (0, 0);
                     \coordinate (2) at (3, 0);
                     \coordinate (3) at (2, 2);
                     \coordinate (4) at (1.5, 0.7);
                     \fill[black!15] (1) -- (4) -- (3) -- cycle;
                     \draw (1) -- (4) -- (3) -- cycle;
                     \draw[dashed] (1) -- (2) -- (3) -- cycle;
                     \draw[dashed] (1) -- (4);
                     \draw[dashed] (2) -- (4);
                     \draw[dashed] (3) -- (4);
                     \foreach \P in {1,3,4} \fill[black] (\P) circle (0.11);
                     \node[left] at (1) {$\mathbf{p}_1$};
                     \node[above] at (3) {$\mathbf{p}_3$};
                     \node[below] at (4) {$\mathbf{p}_4$};
              \end{tikzpicture} \\
              &-
              \begin{tikzpicture}[baseline = 1.0em, x = 1.4em, y = 1.4em]
                     \coordinate (1) at (0, 0);
                     \coordinate (2) at (3, 0);
                     \coordinate (3) at (2, 2);
                     \coordinate (4) at (1.5, 0.7);
                     \draw[dashed] (1) -- (2) -- (3) -- cycle;
                     \draw[thick] (1) -- (4);
                     \draw[dashed] (2) -- (4);
                     \draw[dashed] (3) -- (4);
                     \foreach \P in {1,4} \fill[black] (\P) circle (0.11);
                     \node[left] at (1) {$\mathbf{p}_1$};
                     \node[below] at (4) {$\mathbf{p}_4$};
              \end{tikzpicture} -
              \begin{tikzpicture}[baseline = 1.0em, x = 1.4em, y = 1.4em]
                     \coordinate (1) at (0, 0);
                     \coordinate (2) at (3, 0);
                     \coordinate (3) at (2, 2);
                     \coordinate (4) at (1.5, 0.7);
                     \draw[dashed] (1) -- (2) -- (3) -- cycle;
                     \draw[dashed] (1) -- (4);
                     \draw[thick] (2) -- (4);
                     \draw[dashed] (3) -- (4);
                     \foreach \P in {2,4} \fill[black] (\P) circle (0.11);
                     \node[right] at (2) {$\mathbf{p}_2$};
                     \node[below] at (4) {$\mathbf{p}_4$};
              \end{tikzpicture} -
              \begin{tikzpicture}[baseline = 1.0em, x = 1.4em, y = 1.4em]
                     \coordinate (1) at (0, 0);
                     \coordinate (2) at (3, 0);
                     \coordinate (3) at (2, 2);
                     \coordinate (4) at (1.5, 0.7);
                     \draw[dashed] (1) -- (2) -- (3) -- cycle;
                     \draw[dashed] (1) -- (4);
                     \draw[dashed] (2) -- (4);
                     \draw[thick] (3) -- (4);
                     \foreach \P in {3,4} \fill[black] (\P) circle (0.11);
                     \node[above] at (3) {$\mathbf{p}_3$};
                     \node[below] at (4) {$\mathbf{p}_4$};
              \end{tikzpicture} +
              \begin{tikzpicture}[baseline = 1.0em, x = 1.4em, y = 1.4em]
                     \coordinate (1) at (0, 0);
                     \coordinate (2) at (3, 0);
                     \coordinate (3) at (2, 2);
                     \coordinate (4) at (1.5, 0.7);
                     \draw[dashed] (1) -- (2) -- (3) -- cycle;
                     \draw[dashed] (1) -- (4);
                     \draw[dashed] (2) -- (4);
                     \draw[dashed] (3) -- (4);
                     \foreach \P in {4} \fill[black] (\P) circle (0.11);
                     \node[below] at (4) {$\mathbf{p}_4$};
              \end{tikzpicture}.
       \end{aligned}
       \end{equation}
\end{rem}
We now give the general definition of partitions in $\mathbb{R}^2$. Partitions were introduced by Khovanskii~\cite{Khovanskii1997} and applied by Funke in order to describe partition relations~\cite{Funke2021}. In this paper, we use an extended definition of a partition from that given in~\cite{Khovanskii1997}. The introduction of some notation is necessary in order to define the partitions. Given a certain set of $P = \{\mathbf{p}_1,\dots,\mathbf{p}_n\}$, $\mathcal{P} = \mathrm{conv}(P)$, and $v(P) = \{\mathbf{p}_{v_1},\dots,\mathbf{p}_{v_m}\}$, assume $m \ge 3$ and that the order of indices $v_i$ are given counter-clockwise. Then the boundary of $\mathcal{P}$ is expressed as $\overline{\mathbf{p}_{v_1}\mathbf{p}_{v_2}} \cup \overline{\mathbf{p}_{v_2}\mathbf{p}_{v_3}} \cup \dots \cup \overline{\mathbf{p}_{v_{m-1}}\mathbf{p}_{v_m}} \cup \overline{\mathbf{p}_{v_m}\mathbf{p}_{v_1}}$. We write the set of points of $P$ lying on the boundary of $\mathcal{P}$ as $b(P) := \{\mathbf{p}_{b_1}, \dots, \mathbf{p}_{b_k}\}$, assuming that $\mathbf{p}_{b_1} = \mathbf{p}_{v_1}$ and that the order of indices $b_i$ are given counter-clockwise. Obviously, $v(P) \subset b(P)$. We write $\mathrm{int}(P) := P \setminus b(P)$. We define the set of boundary segments of $\mathcal{P}$ with respect to $P$ as $e(P) := \{\overline{\mathbf{p}_{b_j} \mathbf{p}_{b_{j+1}}} \,|\, j = 1, \dots, k\}$, where $\mathbf{p}_{b_{k+1}}$ is $\mathbf{p}_{b_1}$.

A partition $\Pi := \{P_i \,|\, P_i \subset P, i \in I, |I| < \infty\}$ of $\mathcal{P}$ by $P$ is defined as a finite collection of subsets $P_i$ of $P$ which satisfies the following conditions.
\begin{enumerate}
       \item For any $\mathbf{p} \in P$, there exists at least one $P_i \in \Pi$ such that $\mathbf{p} \in b(P_i)$.
       \item For each $P_i \in \Pi$, $|v(P_i)| \ge 3$.
       \item It holds that $\bigcup_{P_i \in \Pi} \mathrm{conv}(P_i) = \mathrm{conv}(P)$.
       \item If $P_{i_1}, P_{i_2} \in \Pi$ such that $\mathrm{conv}(P_{i_1}) \cap \mathrm{conv}(P_{i_2}) \ne \emptyset$, then $\mathrm{conv}(P_{i_1}) \cap \mathrm{conv}(P_{i_2})$ is a face of $\mathrm{conv}(P_{i_1})$ and $\mathrm{conv}(P_{i_2})$. Note that $\mathcal{Q}$ is a face of $\mathcal{P}$ if $\mathcal{Q} = \mathcal{P}$, $\mathcal{Q} \in e(P)$, or $\mathcal{Q} \in b(P)$.
\end{enumerate}
We define a set of dividing segments of $\Pi$ as $d(\Pi) := (\bigcup_{P_i \in \Pi} e(P_i)) \setminus e(P)$.
\begin{exa}
       In this example, points are indicated by numbers in boldface such as $\mathbf{1}$, $\mathbf{2}$, $\dots$. Let $P = \{\mathbf{1},\mathbf{2},\mathbf{3},\mathbf{4},\mathbf{5},\mathbf{6},\mathbf{7},\mathbf{8},\mathbf{9},\mathbf{10},\mathbf{11}\}$ such that the positions of the points are as shown in Figure \ref{Fig:ExampleOfPartition}. Then $v(P) = \{\mathbf{1},\mathbf{2},\mathbf{3},\mathbf{4},\mathbf{5},\mathbf{6},\mathbf{7}\}$, $b(P) = \{\mathbf{1},\mathbf{2},\mathbf{3},\mathbf{4},\mathbf{5},\mathbf{11},\mathbf{6},\mathbf{7}\}$, and $\mathrm{int}(P) = \{\mathbf{8}, \mathbf{9}, \mathbf{10}\}$. According to our definitions, $\overline{\mathbf{5}\, \mathbf{11}}, \overline{\mathbf{11}\, \mathbf{6}} \in e(P)$, but $\overline{\mathbf{5}\, \mathbf{6}} \notin e(P)$. For this $P$, we can define $\Pi = \{P_1,P_2,P_3,P_4,P_5\}$, where $P_1 = \{\mathbf{1}, \mathbf{8}, \mathbf{7}\}$, $P_2 = \{\mathbf{1}, \mathbf{2}, \mathbf{3}, \mathbf{9}, \mathbf{8}\}$, $P_3 = \{\mathbf{3}, \mathbf{4}, \mathbf{5}, \mathbf{11}, \mathbf{10}, \mathbf{9}\}$, $P_4 = \{\mathbf{6}, \mathbf{7}, \mathbf{8}, \mathbf{10}, \mathbf{11}\}$, $P_5 = \{\mathbf{8}, \mathbf{9}, \mathbf{10}\}$. The dividing segments of $\Pi$ are $\overline{\mathbf{1}\, \mathbf{8}}$, $\overline{\mathbf{3}\, \mathbf{9}}$, $\overline{\mathbf{7}\, \mathbf{8}}$, $\overline{\mathbf{8}\, \mathbf{9}}$, $\overline{\mathbf{8}\, \mathbf{10}}$, $\overline{\mathbf{9}\, \mathbf{10}}$, $\overline{\mathbf{10}\, \mathbf{11}}$.
\end{exa}
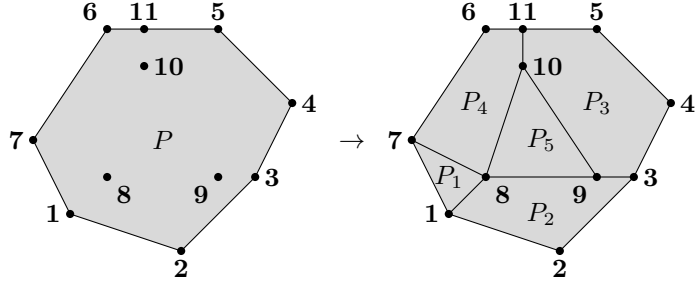
\begin{figure}
       \centering
       \begin{tikzpicture}[baseline = 2.5em, x = 1.4em, y = 1.4em]
              \coordinate[label=left:$\mathbf{1}$] (1) at (0, 0);
              \coordinate[label=below:$\mathbf{2}$] (2) at (3, -1);
              \coordinate[label=right:$\mathbf{3}$] (3) at (5, 1);
              \coordinate[label=right:$\mathbf{4}$] (4) at (6, 3);
              \coordinate[label=above:$\mathbf{5}$] (5) at (4, 5);
              \coordinate[label=above left:$\mathbf{6}$] (6) at (1, 5);
              \coordinate[label=left:$\mathbf{7}$] (7) at (-1, 2);
              \fill[black!15] (1) -- (2) -- (3) -- (4) -- (5) -- (6) -- (7) -- cycle;
              \draw[black] (1) -- (2) -- (3) -- (4) -- (5) -- (6) -- (7) -- cycle;
              \coordinate[label=below right:$\mathbf{8}$] (8) at (1, 1);
              \coordinate[label=below left:$\mathbf{9}$] (9) at (4, 1);
              \coordinate[label=right:$\mathbf{10}$] (10) at (2, 4);
              \coordinate[label=above:$\mathbf{11}$] (11) at (2, 5);
              \foreach \P in {1,2,3,4,5,6,7,8,9,10,11} \fill[black] (\P) circle (0.11);
              \node at (2.5,2) {$P$};
       \end{tikzpicture}
       \,$\rightarrow$\,
       \begin{tikzpicture}[baseline = 2.5em, x = 1.4em, y = 1.4em]
              \coordinate[label=left:$\mathbf{1}$] (1) at (0, 0);
              \coordinate[label=below:$\mathbf{2}$] (2) at (3, -1);
              \coordinate[label=right:$\mathbf{3}$] (3) at (5, 1);
              \coordinate[label=right:$\mathbf{4}$] (4) at (6, 3);
              \coordinate[label=above:$\mathbf{5}$] (5) at (4, 5);
              \coordinate[label=above left:$\mathbf{6}$] (6) at (1, 5);
              \coordinate[label=left:$\mathbf{7}$] (7) at (-1, 2);
              \fill[black!15] (1)--(2)--(3)--(4)--(5)--(6)--(7)--cycle;
              \draw (1)--(2)--(3)--(4)--(5)--(6)--(7)--cycle;
              \coordinate[label=below right:$\mathbf{8}$] (8) at (1, 1);
              \coordinate[label=below left:$\mathbf{9}$] (9) at (4, 1);
              \coordinate[label=right:$\mathbf{10}$] (10) at (2, 4);
              \coordinate[label=above:$\mathbf{11}$] (11) at (2, 5);
              \draw (1)--(8);
              \draw (8)--(9);
              \draw (9)--(10);
              \draw (10)--(8);
              \draw (7)--(8);
              \draw (3)--(9);
              \draw (10)--(11);
              \foreach \P in {1,2,3,4,5,6,7,8,9,10,11} \fill[black] (\P) circle (0.11);
              \node at (0,1) {$P_1$};
              \node at (2.5,0) {$P_2$};
              \node at (4,3) {$P_3$};
              \node at (0.7,3.0) {$P_4$};
              \node at (2.5,2) {$P_5$};
       \end{tikzpicture}
       \caption{An example of a partition of a polygon}
       \label{Fig:ExampleOfPartition}
\end{figure}
\begin{rem}
       Edges of polygons are adopted in the description of partitions in~\cite{Funke2021}. Instead of edges, we adopt boundary segments $e(P)$ which may be a subdivision of edges. Therefore, our definition of a partition is an extension of the original one.
\end{rem}
If $\mathcal{P}$ and its partition $\Pi$ are given, we can derive a signed Minkowski decomposition of $\mathcal{P}$ induced by $\Pi$. It corresponds to the $2$-dimensional case of Funke~\cite[Proposition 3.8]{Funke2021}, and is called the partition relation.
\begin{cor} \label{CorOfGeneralSubdivision}
       Let $P = \{\mathbf{p}_1, \dots, \mathbf{p}_n\}$ and $\Pi$ be a partition of $\mathrm{conv}(P)$. It holds that
       \begin{equation}
              \mathrm{conv}(P) = \sum_{P_i \in \Pi} \mathrm{conv}(P_i) - \sum_{\overline{\mathbf{p}_i \mathbf{p}_j} \in d(\Pi)} \overline{\mathbf{p}_i \mathbf{p}_j} + \sum_{\mathbf{p} \in \mathrm{int}(P)} \mathbf{p}, \label{EqofCorOfGeneralSubdivision}
       \end{equation}
       where $\sum$ represents the Minkowski sum.
\end{cor}
\begin{proof}
       By applying Lemma~\ref{LemOfPolygon} for $\mathrm{conv}(P)$ along each $d(\Pi)$, we obtain, in general, a finer decomposition, such as \eqref{Eq1PrfCorTriangleSubdivision}. Then, by applying Lemma~\ref{LemOfLine} and Lemma~\ref{LemOfPolygon} to this decomposition, we can obtain \eqref{EqofCorOfGeneralSubdivision}. Note that the detailed proof including the case of higher dimensions was already given in~\cite{Funke2021}.
\end{proof}
\begin{rem}
       The following example shows a decomposition which cannot be treated in terms of the original results. Assume $P = \{\mathbf{1} = (1,0), \mathbf{2} = (2,1), \mathbf{3} = (1,2), \mathbf{4} = (0,1), \mathbf{5}= (1,1)\}$ and then there exists a partition $\Pi$ based on our definition given by
       \begin{equation} \label{EqOfNotPolytopalSubdivision}
              \Pi = \{\{\mathbf{1}, \mathbf{2}, \mathbf{5}\},\{\mathbf{1}, \mathbf{4}, \mathbf{5}\},\{\mathbf{2}, \mathbf{3}, \mathbf{4}, \mathbf{5}\}\}
              \Leftrightarrow
              \begin{tikzpicture}[baseline = 1.2em, x = 1.4em, y = 1.4em]
                     \coordinate[label = below:$\mathbf{1}$] (1) at (1, 0);
                     \coordinate[label = right:$\mathbf{2}$] (2) at (2, 1);
                     \coordinate[label = above:$\mathbf{3}$] (3) at (1, 2);
                     \coordinate[label = left:$\mathbf{4}$] (4) at (0, 1);
                     \fill[black!15] (1)--(2)--(3)--(4)--cycle;
                     \coordinate[label = above:$\mathbf{5}$] (5) at (1, 1);
                     \draw (1)--(2)--(3)--(4)--cycle;
                     \draw[dashed] (1)--(5);
                     \draw[dashed] (2)--(4);
                     \foreach \P in {1,2,3,4,5} \fill[black] (\P) circle (0.11);
              \end{tikzpicture}.
       \end{equation}
       This partition \eqref{EqOfNotPolytopalSubdivision} does not satisfy the requirements of the original definition, since $\overline{\mathbf{2}\,\mathbf{5}} = \triangle\mathbf{1}\mathbf{2}\mathbf{5} \cap \triangle\mathbf{2}\mathbf{3}\mathbf{4}$ is not an edge of $\triangle\mathbf{2}\mathbf{3}\mathbf{4}$. While, we can obtain a decomposition from \eqref{EqOfNotPolytopalSubdivision} as follows,
       \begin{equation}
       \begin{tikzpicture}[baseline = 1.2em, x = 1.4em, y = 1.4em]
              \coordinate (1) at (1, 0);
              \coordinate (2) at (2, 1);
              \coordinate (3) at (1, 2);
              \coordinate (4) at (0, 1);
              \coordinate (5) at (1, 1);
              \fill[black!15] (1)--(2)--(3)--(4)--cycle;
              \draw (1)--(2)--(3)--(4)--cycle;
              \draw[dashed] (1)--(5);
              \draw[dashed] (2)--(4);
              \foreach \P in {1,2,3,4,5} \fill[black] (\P) circle (0.11);
       \end{tikzpicture} =
       \begin{tikzpicture}[baseline = 1.2em, x = 1.4em, y = 1.4em]
              \coordinate (1) at (1, 0);
              \coordinate (2) at (2, 1);
              \coordinate (3) at (1, 2);
              \coordinate (4) at (0, 1);
              \coordinate (5) at (1, 1);
              \draw[dashed] (1)--(2)--(3)--(4)--cycle;
              \fill[black!15] (2)--(3)--(4)--cycle;
              \draw (2)--(3)--(4)--cycle;
              \foreach \P in {2,3,4,5} \fill[black] (\P) circle (0.11);
       \end{tikzpicture} +
       \begin{tikzpicture}[baseline = 1.2em, x = 1.4em, y = 1.4em]
              \coordinate (1) at (1, 0);
              \coordinate (2) at (2, 1);
              \coordinate (3) at (1, 2);
              \coordinate (4) at (0, 1);
              \coordinate (5) at (1, 1);
              \fill[black!15] (1)--(5)--(4)--cycle;
              \draw (1)--(5)--(4)--cycle;
              \draw[dashed] (1)--(2)--(3)--(4)--cycle;
              \foreach \P in {1,4,5} \fill[black] (\P) circle (0.11);
       \end{tikzpicture} +
       \begin{tikzpicture}[baseline = 1.2em, x = 1.4em, y = 1.4em]
              \coordinate (1) at (1, 0);
              \coordinate (2) at (2, 1);
              \coordinate (3) at (1, 2);
              \coordinate (4) at (0, 1);
              \coordinate (5) at (1, 1);
              \fill[black!15] (1)--(2)--(5)--cycle;
              \draw (1)--(2)--(5)--cycle;
              \draw[dashed] (1)--(2)--(3)--(4)--cycle;
              \foreach \P in {1,2,5} \fill[black] (\P) circle (0.11);
       \end{tikzpicture} -
       \begin{tikzpicture}[baseline = 1.2em, x = 1.4em, y = 1.4em]
              \coordinate (1) at (1, 0);
              \coordinate (2) at (2, 1);
              \coordinate (3) at (1, 2);
              \coordinate (4) at (0, 1);
              \draw[dashed] (1)--(2)--(3)--(4)--cycle;
              \draw (5)--(4);
              \foreach \P in {4,5} \fill[black] (\P) circle (0.11);
       \end{tikzpicture} -
       \begin{tikzpicture}[baseline = 1.2em, x = 1.4em, y = 1.4em]
              \coordinate (1) at (1, 0);
              \coordinate (2) at (2, 1);
              \coordinate (3) at (1, 2);
              \coordinate (4) at (0, 1);
              \draw[dashed] (1)--(2)--(3)--(4)--cycle;
              \draw (1)--(5);
              \foreach \P in {1,5} \fill[black] (\P) circle (0.11);
       \end{tikzpicture} -
       \begin{tikzpicture}[baseline = 1.2em, x = 1.4em, y = 1.4em]
              \coordinate (1) at (1, 0);
              \coordinate (2) at (2, 1);
              \coordinate (3) at (1, 2);
              \coordinate (4) at (0, 1);
              \draw[dashed] (1)--(2)--(3)--(4)--cycle;
              \draw (2)--(5);
              \foreach \P in {2,5} \fill[black] (\P) circle (0.11);
       \end{tikzpicture} +
       \begin{tikzpicture}[baseline = 1.2em, x = 1.4em, y = 1.4em]
              \coordinate (1) at (1, 0);
              \coordinate (2) at (2, 1);
              \coordinate (3) at (1, 2);
              \coordinate (4) at (0, 1);
              \draw[dashed] (1)--(2)--(3)--(4)--cycle;
              \foreach \P in {5} \fill[black] (\P) circle (0.11);
       \end{tikzpicture},
       \end{equation}
       by applying the cutting relations. Though the geometric condition is extended, the proof of Corollary~\ref{CorOfGeneralSubdivision} follows from the original one~\cite{Funke2021}.
\end{rem}
\section{Decompositions of integral polygons} \label{Sec:DecompOfLatPolytopes}
Consider a polygon $\mathcal{P}$ with vertices $v(\mathcal{P})$. If all coordinates of $\mathbf{p}_{v_i} \in v(\mathcal{P})$ are integers, $\mathcal{P}$ is called an integral polygon. We will show that any integral polygon can be represented as a signed Minkowski decomposition in a specific form.

Consider integral triangles. The minimum value of area of an integral triangle is $1/2$, and we call such a triangle a minimum triangle. Let $T_{\mathrm{uni}}$ denote the set of all integral triangles $\mathcal{T} = \mathrm{conv}\{(a_1,b_1), (a_2,b_2), (a_3,b_3)\}$ which satisfy the following conditions.
\begin{enumerate}
       \item (Nonnegativity) The coordinates of the three vertices $(a_1,b_1)$, $(a_2,b_2)$, $(a_3,b_3)$ are all non-negative.
       \item (Minimality) The area of $\mathcal{T}$ is $1/2$, that is, it is a minimum triangle.
       \item At least one of $a_1, a_2, a_3$ is $0$, and at least one of $b_1, b_2, b_3$ is $0$.
\end{enumerate}       
In other words, $\mathcal{T} \in T_{\mathrm{uni}}$ is given in the form of $\mathrm{conv}\{(0,0), (*,*), (*,*)\}$, or $\mathrm{conv}\{(*,0), (0,*), (*,*)\}$, where the values of $*$ are taken from nonnegative integers as long as $\mathcal{T}$ is a minimum triangle. We call an element of $T_{\mathrm{uni}}$ a unit triangle. Note that any minimum triangle can be uniquely represented as a translated unit triangle.
\begin{rem}
       Any minimum triangle does not contain any integral points in the boundary or interior, except at the vertices.
\end{rem}
We define
\begin{align}
       \mathcal{I}_x := \mathrm{conv}\{(0,0),(1,0)\} =
       \begin{tikzpicture}[baseline = 0.6em, x = 1.4em, y = 1.4em]
              \coordinate (1) at (0, 0);
              \coordinate (2) at (1, 0);
              \draw[step = 1, black, dashed] (-0.2,-0.2) grid (1.2,1.2);
              \draw[thick] (1) -- (2);
              \foreach \P in {1,2} \fill[black] (\P) circle (0.11);
       \end{tikzpicture}, \\
       \mathcal{I}_y := \mathrm{conv}\{(0,0),(0,1)\} =
       \begin{tikzpicture}[baseline = 0.6em, x = 1.4em, y = 1.4em]
              \coordinate (1) at (0, 0);
              \coordinate (2) at (0, 1);
              \draw[step = 1, black, dashed] (-0.2,-0.2) grid (1.2,1.2);
              \draw[thick] (1) -- (2);
              \foreach \P in {1,2} \fill[black] (\P) circle (0.11);
       \end{tikzpicture},
\end{align}
and call $\mathcal{I}_x$ a $x$-unit segment and $\mathcal{I}_y$ a $y$-unit segment respectively.

Any integral polygon can be represented as a signed decomposition, consisting of $\mathcal{I}_x$, $\mathcal{I}_y$, and unit triangles.
\begin{prp} \label{PropOfDecompFormula}
       For any integral segment or polygon $\mathcal{P}$, there exists a decomposition of $\mathcal{P}$ expressed as
       \begin{equation} \label{DecompFormula}
              \mathcal{P} = \mathbf{t} + k_x \mathcal{I}_x + k_y \mathcal{I}_y +\sum_{\mathcal{T}_i \in S} k_{\mathcal{T}_i} \mathcal{T}_i,
       \end{equation}
       where $\mathbf{t} \in \mathbb{Z}^2$, $k_x, k_y, k_{\mathcal{T}_i}$ are integers, and $S$ is a finite subset of $T_{\mathrm{uni}}$.
\end{prp}
Before the proof of Proposition~\ref{PropOfDecompFormula}, we will show some example decompositions.
\begin{exa} \label{ExamplesOfUnitDecomp}
       We can obtain the following decompositions.
       \begin{equation}
              \begin{tikzpicture}[baseline = 0.6em, x = 1.4em, y = 1.4em]
                     \coordinate (1) at (0, 0);
                     \coordinate (2) at (2, 0);
                     \coordinate (3) at (2, 2);
                     \draw[thick] (1) -- (2) -- (3) -- cycle;
                     \fill[black!15] (1) -- (2) -- (3) -- cycle;
                     \draw[step = 1, black, dashed] (-0.2,-0.2) grid (2.2,2.2);
                     \foreach \P in {1,2,3} \fill[black] (\P) circle (0.11);
              \end{tikzpicture}
              = 2
              \begin{tikzpicture}[baseline = 0.6em, x = 1.4em, y = 1.4em]
                     \coordinate (1) at (0, 0);
                     \coordinate (2) at (1, 0);
                     \coordinate (3) at (1, 1);
                     \draw[thick] (1) -- (2) -- (3) -- cycle;
                     \fill[black!15] (1) -- (2) -- (3) -- cycle;
                     \draw[step = 1, black, dashed] (-0.2,-0.2) grid (1.2,1.2);
                     \foreach \P in {1,2,3} \fill[black] (\P) circle (0.11);
              \end{tikzpicture},
       \end{equation}
       \begin{equation}
              \begin{tikzpicture}[baseline = -0.8em, x = 1.4em, y = 1.4em]
                     \coordinate (1) at (-1, -1);
                     \coordinate (2) at (1, 0);
                     \draw[step = 1, black, dashed] (-1.2,-1.2) grid (1.2,0.2);
                     \draw[thick] (1) -- (2);
                     \foreach \P in {1,2} \fill[black] (\P) circle (0.11);
              \end{tikzpicture} = 
              \begin{tikzpicture}[baseline = 0.6em, x = 1.4em, y = 1.4em]
                     \coordinate (1) at (0, 0);
                     \coordinate (2) at (0, 1);
                     \draw[step = 1, black, dashed] (-0.2,-0.2) grid (1.2,1.2);
                     \draw[thick] (1) -- (2);
                     \foreach \P in {1,2} \fill[black] (\P) circle (0.11);
              \end{tikzpicture} -
              \begin{tikzpicture}[baseline = 0.6em, x = 1.4em, y = 1.4em]
                     \coordinate (1) at (0, 0);
                     \coordinate (2) at (1, 0);
                     \coordinate (3) at (1, 1);
                     \draw[thick] (1) -- (2) -- (3) -- cycle;
                     \fill[black!15] (1) -- (2) -- (3) -- cycle;
                     \draw[step = 1, black, dashed] (-0.2,-0.2) grid (1.2,1.2);
                     \foreach \P in {1,2,3} \fill[black] (\P) circle (0.11);
              \end{tikzpicture} -
              \begin{tikzpicture}[baseline = 0.6em, x = 1.4em, y = 1.4em]
                     \coordinate (1) at (0, 0);
                     \coordinate (2) at (0, 1);
                     \coordinate (3) at (1, 1);
                     \draw[thick] (1) -- (2) -- (3) -- cycle;
                     \fill[black!15] (1) -- (2) -- (3) -- cycle;
                     \draw[step = 1, black, dashed] (-0.2,-0.2) grid (1.2,1.2);
                     \foreach \P in {1,2,3} \fill[black] (\P) circle (0.11);
              \end{tikzpicture} +
              \begin{tikzpicture}[baseline = 0.6em, x = 1.4em, y = 1.4em]
                     \coordinate (1) at (0, 0);
                     \coordinate (2) at (1, 0);
                     \coordinate (3) at (2, 1);
                     \draw[thick] (1) -- (2) -- (3) -- cycle;
                     \fill[black!15] (1) -- (2) -- (3) -- cycle;
                     \draw[step = 1, black, dashed] (-0.2,-0.2) grid (2.2,1.2);
                     \foreach \P in {1,2,3} \fill[black] (\P) circle (0.11);
              \end{tikzpicture} +
              \begin{tikzpicture}[baseline = 0.6em, x = 1.4em, y = 1.4em]
                     \coordinate (1) at (0, 0);
                     \coordinate (2) at (1, 1);
                     \coordinate (3) at (2, 1);
                     \draw[thick] (1) -- (2) -- (3) -- cycle;
                     \fill[black!15] (1) -- (2) -- (3) -- cycle;
                     \draw[step = 1, black, dashed] (-0.2,-0.2) grid (2.2,1.2);
                     \foreach \P in {1,2,3} \fill[black] (\P) circle (0.11);
              \end{tikzpicture},
       \end{equation}
       \begin{equation}
              \label{EqOfExamplesOfUnitDecomp}
              \begin{tikzpicture}[baseline = 0.6em, x = 1.4em, y = 1.4em]
                     \coordinate (1) at (0, 0);
                     \coordinate (3) at (1, 1);
                     \coordinate (4) at (3, 0);
                     \fill[black!15] (1) -- (3) -- (4) -- cycle;
                     \draw[step = 1, black, dashed] (-0.2,-0.2) grid (3.2,1.2);
                     \draw (1) -- (3) -- (4) -- cycle;
                     \foreach \P in {1,3,4} \fill[black] (\P) circle (0.11);
              \end{tikzpicture} =
              \begin{tikzpicture}[baseline = 0.6em, x = 1.4em, y = 1.4em]
                     \coordinate (1) at (0, 0);
                     \coordinate (2) at (1, 0);
                     \draw[step = 1, black, dashed] (-0.2,-0.2) grid (1.2,1.2);
                     \draw[thick] (1) -- (2);
                     \foreach \P in {1,2} \fill[black] (\P) circle (0.11);
              \end{tikzpicture} +
              \begin{tikzpicture}[baseline = 0.6em, x = 1.4em, y = 1.4em]
                     \coordinate (1) at (0, 0);
                     \coordinate (2) at (1, 0);
                     \coordinate (3) at (1, 1);
                     \fill[black!15] (1) -- (2) -- (3) -- cycle;
                     \draw[step = 1, black, dashed] (-0.2,-0.2) grid (1.2,1.2);
                     \draw (1) -- (2) -- (3) -- cycle;
                     \foreach \P in {1,2,3} \fill[black] (\P) circle (0.11);
              \end{tikzpicture} +
              \begin{tikzpicture}[baseline = 0.6em, x = 1.4em, y = 1.4em]
                     \coordinate (1) at (1, 0);
                     \coordinate (2) at (0, 1);
                     \coordinate (3) at (2, 0);
                     \fill[black!15] (1) -- (2) -- (3) -- cycle;
                     \draw[step = 1, black, dashed] (-0.2,-0.2) grid (2.2,1.2);
                     \draw (1) -- (2) -- (3) -- cycle;
                     \foreach \P in {1,2,3} \fill[black] (\P) circle (0.11);
              \end{tikzpicture} -
              \begin{tikzpicture}[baseline = 0.6em, x = 1.4em, y = 1.4em]
                     \coordinate (1) at (1, 0);
                     \coordinate (2) at (0, 1);
                     \coordinate (3) at (1, 1);
                     \fill[black!15] (1) -- (2) -- (3) -- cycle;
                     \draw[step = 1, black, dashed] (-0.2,-0.2) grid (1.2,1.2);
                     \draw (1) -- (2) -- (3) -- cycle;
                     \foreach \P in {1,2,3} \fill[black] (\P) circle (0.11);
              \end{tikzpicture},
       \end{equation}
       \begin{equation}
              \begin{tikzpicture}[baseline = 0.6em, x = 1.4em, y = 1.4em]
                     \coordinate (1) at (1, 0);
                     \coordinate (2) at (0, 1);
                     \coordinate (3) at (1, 3);
                     \coordinate (4) at (4, 1);
                     \fill[black!15] (1) -- (2) -- (3) -- (4) -- cycle;
                     \draw[step = 1, black, dashed] (-0.2,-0.2) grid (4.2,3.2);
                     \draw (1) -- (2) -- (3) -- (4) -- cycle;
                     \foreach \P in {1,2,3,4} \fill[black] (\P) circle (0.11);
              \end{tikzpicture} = -
              \begin{tikzpicture}[baseline = 0.6em, x = 1.4em, y = 1.4em]
                     \coordinate (1) at (0, 0);
                     \coordinate (2) at (0, 1);
                     \draw[step = 1, black, dashed] (-0.2,-0.2) grid (1.2,1.2);
                     \draw[thick] (1) -- (2);
                     \foreach \P in {1,2} \fill[black] (\P) circle (0.11);
              \end{tikzpicture} +
              \begin{tikzpicture}[baseline = 0.6em, x = 1.4em, y = 1.4em]
                     \coordinate (1) at (0, 0);
                     \coordinate (2) at (1, 0);
                     \coordinate (3) at (0, 1);
                     \fill[black!15] (1) -- (2) -- (3) -- cycle;
                     \draw[step = 1, black, dashed] (-0.2,-0.2) grid (1.2,1.2);
                     \draw (1) -- (2) -- (3) -- cycle;
                     \foreach \P in {1,2,3} \fill[black] (\P) circle (0.11);
              \end{tikzpicture} -
              \begin{tikzpicture}[baseline = 0.6em, x = 1.4em, y = 1.4em]
                     \coordinate (1) at (0, 0);
                     \coordinate (2) at (1, 0);
                     \coordinate (3) at (2, 1);
                     \fill[black!15] (1) -- (2) -- (3) -- cycle;
                     \draw[step = 1, black, dashed] (-0.2,-0.2) grid (2.2,1.2);
                     \draw (1) -- (2) -- (3) -- cycle;
                     \foreach \P in {1,2,3} \fill[black] (\P) circle (0.11);
              \end{tikzpicture} -
              \begin{tikzpicture}[baseline = 0.6em, x = 1.4em, y = 1.4em]
                     \coordinate (1) at (0, 1);
                     \coordinate (2) at (1, 1);
                     \coordinate (3) at (2, 0);
                     \fill[black!15] (1) -- (2) -- (3) -- cycle;
                     \draw[step = 1, black, dashed] (-0.2,-0.2) grid (2.2,1.2);
                     \draw (1) -- (2) -- (3) -- cycle;
                     \foreach \P in {1,2,3} \fill[black] (\P) circle (0.11);
              \end{tikzpicture} +
              \begin{tikzpicture}[baseline = 0.6em, x = 1.4em, y = 1.4em]
                     \coordinate (1) at (0, 0);
                     \coordinate (2) at (2, 1);
                     \coordinate (3) at (3, 1);
                     \fill[black!15] (1) -- (2) -- (3) -- cycle;
                     \draw[step = 1, black, dashed] (-0.2,-0.2) grid (3.2,1.2);
                     \draw (1) -- (2) -- (3) -- cycle;
                     \foreach \P in {1,2,3} \fill[black] (\P) circle (0.11);
              \end{tikzpicture} +
              \begin{tikzpicture}[baseline = 0.6em, x = 1.4em, y = 1.4em]
                     \coordinate (1) at (0, 0);
                     \coordinate (2) at (1, 1);
                     \coordinate (3) at (1, 2);
                     \fill[black!15] (1) -- (2) -- (3) -- cycle;
                     \draw[step = 1, black, dashed] (-0.2,-0.2) grid (1.2,2.2);
                     \draw (1) -- (2) -- (3) -- cycle;
                     \foreach \P in {1,2,3} \fill[black] (\P) circle (0.11);
              \end{tikzpicture} +
              \begin{tikzpicture}[baseline = 0.6em, x = 1.4em, y = 1.4em]
                     \coordinate (1) at (1, 1);
                     \coordinate (2) at (0, 2);
                     \coordinate (3) at (3, 0);
                     \fill[black!15] (1) -- (2) -- (3) -- cycle;
                     \draw[step = 1, black, dashed] (-0.2,-0.2) grid (3.2,2.2);
                     \draw (1) -- (2) -- (3) -- cycle;
                     \foreach \P in {1,2,3} \fill[black] (\P) circle (0.11);
              \end{tikzpicture},
       \end{equation}
\end{exa}
\noindent All of these decompositions can be verified by straightforward calculations.
\subsection{The proof of Proposition~\ref{PropOfDecompFormula}}
First, we will prove Proposition~\ref{PropOfDecompFormula} for the case of segments. Suppose $\mathbf{p}, \mathbf{q} \in \mathbb{Z}^2$ and $\mathbf{p} \ne \mathbf{q}$. If $\overline{\mathbf{p}\mathbf{q}}$ does not contain any other integral points, we call it a prime segment, otherwise non-prime. According to Lemma~\ref{LemOfLine}, non-prime segments can be decomposed into prime segments and points. Therefore, it is sufficient to consider just prime segments. Additionally, we also assume without loss of generality that one of endpoints of the considered segment is fixed at the origin $\mathbf{o} := (0,0)$, otherwise the considered segment can be translated by an integral vector.
\begin{lem} \label{LemOfEdgesEqualToTriangles}
       Given $\mathbf{p}, \mathbf{q} \in \mathbb{R}^2$, and $\mathbf{r} := \mathbf{p} - \mathbf{q}$. It holds that
       \begin{equation} \label{EquiOfEdgesEqualToTriangles}
              \begin{aligned}
                     \overline{\mathbf{o}\mathbf{p}} + \overline{\mathbf{o}\mathbf{q}} + \overline{\mathbf{o}\mathbf{r}} =& \mathrm{conv}\{\mathbf{o},\mathbf{p},\mathbf{q}\} + \mathrm{conv}\{\mathbf{o},\mathbf{p},\mathbf{r}\} \\
                     =& \mathrm{conv}\{\mathbf{o},\mathbf{q},\mathbf{p}+\mathbf{q},\mathbf{p}+\mathbf{p},\mathbf{p}+\mathbf{r},\mathbf{r}\}.
              \end{aligned}
       \end{equation}
\end{lem}
\begin{proof}
       Equation \eqref{EquiOfEdgesEqualToTriangles} is derived from the following relation, including the case where $\overline{\mathbf{o} \mathbf{p}}$ is parallel with $\overline{\mathbf{o} \mathbf{q}}$,
       \begin{equation} \label{EquiOfPolygonOfEdges}
              \begin{tikzpicture}[baseline = 0.5em, x = 1.4em, y = 1.4em]
                     \coordinate (1) at (0, 0);
                     \node[below] at (1) {$\mathbf{o}$};
                     \coordinate (2) at (2, 1);
                     \node[above] at (2) {$\mathbf{p}$};
                     \coordinate (3) at (1, 0);
                     \coordinate (4) at (1, 1);
                     \draw[dashed] (1) -- (2) -- (3) -- cycle;
                     \draw[dashed] (1) -- (2) -- (4) -- cycle;
                     \draw[thick] (1) -- (2);
                     \foreach \P in {1,2} \fill[black] (\P) circle (0.11);
              \end{tikzpicture} +
              \begin{tikzpicture}[baseline = 0.5em, x = 1.4em, y = 1.4em]
                     \coordinate (1) at (0, 0);
                     \node[below] at (1) {$\mathbf{o}$};
                     \coordinate (2) at (2, 1);
                     \coordinate (3) at (1, 0);
                     \node[below] at (3) {$\mathbf{q}$};
                     \coordinate (4) at (1, 1);
                     \draw[dashed] (1) -- (2) -- (3) -- cycle;
                     \draw[dashed] (1) -- (2) -- (4) -- cycle;
                     \draw[thick] (1) -- (3);
                     \foreach \P in {1,3} \fill[black] (\P) circle (0.11);
              \end{tikzpicture} +
              \begin{tikzpicture}[baseline = 0.5em, x = 1.4em, y = 1.4em]
                     \coordinate (1) at (0, 0);
                     \node[below] at (1) {$\mathbf{o}$};
                     \coordinate (2) at (2, 1);
                     \coordinate (3) at (1, 0);
                     \coordinate (4) at (1, 1);
                     \node[above] at (4) {$\mathbf{r}$};
                     \draw[dashed] (1) -- (2) -- (3) -- cycle;
                     \draw[dashed] (1) -- (2) -- (4) -- cycle;
                     \draw[thick] (1) -- (4);
                     \foreach \P in {1,4} \fill[black] (\P) circle (0.11);
              \end{tikzpicture} =
              \begin{tikzpicture}[baseline = 0.5em, x = 1.4em, y = 1.4em]
                     \coordinate (1) at (0, 0);
                     \node[below] at (1) {$\mathbf{o}$};
                     \coordinate (2) at (2, 1);
                     \node[above] at (2) {$\mathbf{p}$};
                     \coordinate (3) at (1, 0);
                     \node[below] at (3) {$\mathbf{q}$};
                     \coordinate (4) at (1, 1);
                     \fill[black!15] (1) -- (2) -- (3) -- cycle;
                     \draw[dashed] (1) -- (2) -- (4) -- cycle;
                     \draw[thick] (1) -- (2) -- (3) -- cycle;
                     \foreach \P in {1,2,3} \fill[black] (\P) circle (0.11);
              \end{tikzpicture} +
              \begin{tikzpicture}[baseline = 0.5em, x = 1.4em, y = 1.4em]
                     \coordinate (1) at (0, 0);
                     \node[below] at (1) {$\mathbf{o}$};
                     \coordinate (2) at (2, 1);
                     \node[above] at (2) {$\mathbf{p}$};
                     \coordinate (3) at (1, 0);
                     \coordinate (4) at (1, 1);
                     \node[above] at (4) {$\mathbf{r}$};
                     \fill[black!15] (1) -- (2) -- (4) -- cycle;
                     \draw[dashed] (1) -- (2) -- (3) -- cycle;
                     \draw[thick] (1) -- (2) -- (4) -- cycle;
                     \foreach \P in {1,2,4} \fill[black] (\P) circle (0.11);
              \end{tikzpicture}
              =
              \begin{tikzpicture}[baseline = 1.6em, x = 1.4em, y = 1.4em]
                     \coordinate (1) at (0, 0);
                     \coordinate (2) at (1, 0);
                     \coordinate (3) at (3, 1);
                     \coordinate (4) at (4, 2);
                     \coordinate (5) at (3, 2);
                     \coordinate (6) at (1, 1);
                     \coordinate (7) at (2, 1);
                     \node[below left] at (1) {$\mathbf{o}$};
                     \node[below] at (2) {$\mathbf{q}$};
                     \node[right] at (3) {$\mathbf{p}+\mathbf{q}$};
                     \node[right] at (4) {$\mathbf{p} + \mathbf{p}$};
                     \node[above] at (5) {$\mathbf{p}+\mathbf{r}$};
                     \node[left] at (6) {$\mathbf{r}$};
                     \fill[black!15] (1) -- (2) -- (3) -- (4) -- (5) -- (6) -- cycle;
                     \draw[thick] (1) -- (2) -- (3) -- (4) -- (5) -- (6) -- cycle;
                     \foreach \P in {1,2,3,4,5,6,7} \fill[black] (\P) circle (0.11);
                     \draw[dashed] (1) -- (4);
                     \draw[dashed] (2) -- (5);
                     \draw[dashed] (3) -- (6);
              \end{tikzpicture}.
       \end{equation}
\end{proof}
\begin{lem} \label{LemOfShorterExist}
       Assume that $\overline{\mathbf{o}\mathbf{p}}$ is a prime segment, and $\mathbf{p} := (a,b)$ is not $(\pm1,0)$ or $(0,\pm1)$. Then there exists a point $\mathbf{q} \in \mathbb{Z}^2$ such that
       \begin{enumerate}
              \item The segment $\overline{\mathbf{o}\mathbf{q}}$ is prime, and its length is shorter than that of $\overline{\mathbf{o}\mathbf{p}}$. \label{LengthCond}
              \item The triangle $\triangle \mathbf{o} \mathbf{p} \mathbf{q}$ is minimum. \label{AreaCond}
       \end{enumerate}
       Additionally, if $\mathbf{r} := \mathbf{p} - \mathbf{q}$, then it holds that
       \begin{enumerate}
              \item The segment $\overline{\mathbf{o}\mathbf{r}}$ is prime, and its length is shorter than that of $\overline{\mathbf{o}\mathbf{p}}$.
              \item The triangle $\triangle \mathbf{o} \mathbf{p} \mathbf{r}$ is minimum.
       \end{enumerate}
\end{lem}
\begin{proof}
       Assume $a > 0$ and $b > 0$ without loss of generality. Consider the set of integral points $I = \{(i,j) \,;\, i,j \in \mathbb{Z},\, 0 \le i \le a,\, 0 \le j \le b\}$. This set is obviously finite. For any $\mathbf{i} \in I$, the length of $\overline{\mathbf{o} \mathbf{i}}$ is shorter than that of $\overline{\mathbf{o} \mathbf{p}}$. There must exist at least one point $(a',b') \in I$ such that the triangle $\mathrm{conv}\{\mathbf{o},\mathbf{p}, (a', b')\}$ does not contain any integral points except its vertices, $\mathbf{o}$, $\mathbf{p}$, $(a',b')$. Therefore, we can choose $(a', b')$ as $\mathbf{q}$. As a result of the symmetry, $\mathbf{r} = \mathbf{p} - \mathbf{q} \in I$ satisfies the same properties as $\mathbf{q}$.
\end{proof}
\begin{figure}
       \centering
       \begin{tikzpicture}[baseline = 0.6em, x = 2.4em, y = 2.4em]
              \draw[->,thick] (-0.4,0) -- (4.4,0);
              \draw[->,thick] (0,-0.4) -- (0,3.4);
              \draw[dashed] (0,3) -- (4,3);
              \draw[dashed] (4,0) -- (4,3);
              \draw (0,0) -- (3,2);
              \draw (4,3) -- (3,2);                  
              \draw (0,0) -- (1,1);                  
              \draw (4,3) -- (1,1);                                    
              \draw[thick] (0,0) -- (4,3);
              \foreach \P in {0,1,2,3,4} \fill[black] (\P,0) circle (0.11);
              \foreach \P in {0,1,2,3,4} \fill[black] (\P,1) circle (0.11);
              \foreach \P in {0,1,2,3,4} \fill[black] (\P,2) circle (0.11);
              \foreach \P in {0,1,2,3,4} \fill[black] (\P,3) circle (0.11);
              \fill[white] (3,2) circle (0.08);
              \fill[white] (1,1) circle (0.08);
              \node[above right] at (4,3) {$\mathbf{p}$};
              \node[below] at (3,2) {$\mathbf{q}$};
              \node[above] at (1,1) {$\mathbf{p}-\mathbf{q}$};
              \node[below] at (4,0) {$a$};
              \node[left] at (0,3) {$b$};
              \node[below left] at (0,0) {$\mathbf{o}$};
       \end{tikzpicture}
       \caption{Lattice points of $I$ and the triangles}
       \label{FigOfPointsInCircle}
\end{figure}
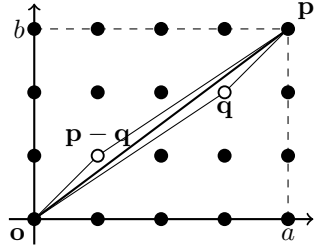
\begin{cor} \label{CorOfEdgeToShorterAndTriangles}
       If $\overline{\mathbf{o}\mathbf{p}}$ is a prime segment, where $\mathbf{p} := (a,b) \ne (\pm1,0)$ or $(0,\pm1)$, then $\overline{\mathbf{o}\mathbf{p}}$ can be decomposed into
       \begin{equation} \label{DecompOf0ax+by}
       \begin{gathered}
              \overline{\mathbf{o}\mathbf{p}} = \triangle \mathbf{o} \mathbf{p} \mathbf{q} + \triangle \mathbf{o} \mathbf{p} \mathbf{r} - \overline{\mathbf{o}\mathbf{q}} - \overline{\mathbf{o}\mathbf{r}}, \\
              \Leftrightarrow
              \begin{tikzpicture}[baseline = 0.6em, x = 1.4em, y = 1.4em]
                     \coordinate (1) at (0, 0);
                     \coordinate (2) at (2, 1);
                     \node[below] at (1) {$\mathbf{o}$};
                     \node[above] at (2) {$\mathbf{p}$};
                     \coordinate (3) at (1, 0);
                     \coordinate (4) at (1, 1);
                     \draw[dashed] (1) -- (2) -- (3) -- cycle;
                     \draw[dashed] (1) -- (2) -- (4) -- cycle;
                     \draw[thick] (1) -- (2);
                     \foreach \P in {1,2} \fill[black] (\P) circle (0.11);
              \end{tikzpicture} =
              \begin{tikzpicture}[baseline = 0.6em, x = 1.4em, y = 1.4em]
                     \coordinate (1) at (0, 0);
                     \coordinate (2) at (2, 1);
                     \coordinate (3) at (1, 0);
                     \coordinate (4) at (1, 1);
                     \fill[black!15] (1) -- (2) -- (3) -- cycle;
                     \draw[dashed] (1) -- (2) -- (4) -- cycle;
                     \draw[thick] (1) -- (2) -- (3) -- cycle;
                     \foreach \P in {1,2,3} \fill[black] (\P) circle (0.11);
                     \node[below] at (1) {$\mathbf{o}$};
                     \node[above] at (2) {$\mathbf{p}$};
                     \node[below] at (3) {$\mathbf{q}$};
              \end{tikzpicture} +
              \begin{tikzpicture}[baseline = 0.6em, x = 1.4em, y = 1.4em]
                     \coordinate (1) at (0, 0);
                     \coordinate (2) at (2, 1);
                     \coordinate (3) at (1, 0);
                     \coordinate (4) at (1, 1);
                     \fill[black!15] (1) -- (2) -- (4) -- cycle;
                     \draw[dashed] (1) -- (2) -- (3) -- cycle;
                     \draw[thick] (1) -- (2) -- (4) -- cycle;
                     \foreach \P in {1,2,4} \fill[black] (\P) circle (0.11);
                     \node[below] at (1) {$\mathbf{o}$};
                     \node[above] at (2) {$\mathbf{p}$};
                     \node[above] at (4) {$\mathbf{r}$};
              \end{tikzpicture} -
              \begin{tikzpicture}[baseline = 0.6em, x = 1.4em, y = 1.4em]
                     \coordinate (1) at (0, 0);
                     \coordinate (2) at (2, 1);
                     \coordinate (3) at (1, 0);
                     \coordinate (4) at (1, 1);
                     \draw[dashed] (1) -- (2) -- (3) -- cycle;
                     \draw[dashed] (1) -- (2) -- (4) -- cycle;
                     \draw[thick] (1) -- (3);
                     \foreach \P in {1,3} \fill[black] (\P) circle (0.11);
                     \node[below] at (1) {$\mathbf{o}$};
                     \node[below] at (3) {$\mathbf{q}$};
              \end{tikzpicture} -
              \begin{tikzpicture}[baseline = 0.6em, x = 1.4em, y = 1.4em]
                     \coordinate (1) at (0, 0);
                     \coordinate (2) at (2, 1);
                     \coordinate (3) at (1, 0);
                     \coordinate (4) at (1, 1);
                     \draw[dashed] (1) -- (2) -- (3) -- cycle;
                     \draw[dashed] (1) -- (2) -- (4) -- cycle;
                     \draw[thick] (1) -- (4);
                     \foreach \P in {1,4} \fill[black] (\P) circle (0.11);
                     \node[below] at (1) {$\mathbf{o}$};
                     \node[above] at (4) {$\mathbf{r}$};
              \end{tikzpicture},
       \end{gathered}
       \end{equation}
       where $\mathbf{q} := (a',b')$, $\mathbf{r} := \mathbf{p} - \mathbf{q}$ such that $\triangle \mathbf{o} \mathbf{p} \mathbf{q}$ and $\triangle \mathbf{o} \mathbf{p} \mathbf{r}$ are minimum, $\overline{\mathbf{o}\mathbf{q}}$ and $\overline{\mathbf{o}\mathbf{r}}$ are prime, and their lengths are shorter than that of $\overline{\mathbf{o}\mathbf{p}}$.
\end{cor}
\begin{proof}
       It follows immediately from Lemma~\ref{LemOfEdgesEqualToTriangles} and Lemma~\ref{LemOfShorterExist}.
\end{proof}
By applying Corollary~\ref{CorOfEdgeToShorterAndTriangles} recursively, Proposition~\ref{PropOfDecompFormula} for segments can be proven as follows.
\begin{lem} \label{LemOfSegmentDecomposition}
       For any segment $\overline{\mathbf{p}\mathbf{q}}$, there exists a decomposition,
       \begin{equation}
              \overline{\mathbf{p}\mathbf{q}} = \mathbf{t} + k_x \mathcal{I}_x + k_y \mathcal{I}_y +\sum_{\mathcal{T} \in S} k_{\mathcal{T}} \mathcal{T},
       \end{equation}
       where the symbols are as given in Proposition \ref{PropOfDecompFormula}.
\end{lem}
\begin{proof}
       Assume that $\overline{\mathbf{p}\mathbf{q}}$ is not prime and is represented as 
       \begin{equation}
              \begin{aligned}
                     \overline{\mathbf{p}\mathbf{q}} = \mathrm{conv}\{&\mathbf{p}_0 := \mathbf{p}, \mathbf{p}_1 := \mathbf{p}_0 + \frac{1}{n} (\mathbf{q} - \mathbf{p}), \mathbf{p}_2 := \mathbf{p}_0 + \frac{2}{n} (\mathbf{q} - \mathbf{p}), \\
                     &\dots, \mathbf{p}_{n-1} := \mathbf{p}_0 + \frac{n-1}{n} (\mathbf{q} - \mathbf{p}), \mathbf{p}_n := \mathbf{q}\},
              \end{aligned}
       \end{equation}
       where $2 \le n \in \mathbb{Z}$, and $\mathbf{p}_1$ to $\mathbf{p}_{n-1}$ are integral points on $\overline{\mathbf{p}\mathbf{q}}$. According to the cutting relation,
       \begin{equation}
              \mathrm{conv}\{\mathbf{p}_0, \mathbf{p}_1, \mathbf{p}_2, \dots, \mathbf{p}_n\} = \sum_{i = 0}^{n-1} \overline{\mathbf{p}_i \mathbf{p}_{i+1}} - \sum_{i = 1}^{n-1} \mathbf{p}_i.
       \end{equation}
       By translation, each $\overline{\mathbf{p}_i \mathbf{p}_{i+1}}$ can be rewritten as $\overline{\mathbf{p}_i \mathbf{p}_{i+1}} = \overline{\mathbf{o} \mathbf{p}_{i}'} + \mathbf{p}_i$, where $\mathbf{p}_i' := \mathbf{p}_{i+1} - \mathbf{p}_i$. Then
       \begin{equation}
              \begin{aligned}
                     \sum_{i = 0}^{n-1} \overline{\mathbf{p}_i \mathbf{p}_{i+1}} - \sum_{i = 1}^{n-1} \mathbf{p}_i
                     &= \sum_{i = 0}^{n-1} \overline{\mathbf{o} \mathbf{p}_{i}'} + \sum_{i = 0}^{n-1} \mathbf{p}_i - \sum_{i = 1}^{n-1} \mathbf{p}_i \\
                     &= \sum_{i = 0}^{n-1} \overline{\mathbf{o} \mathbf{p}_{i}'} + \mathbf{p}_0
              \end{aligned}
       \end{equation}
       is obtained. Since $\overline{\mathbf{o} \mathbf{p}_{0}'} = \overline{\mathbf{o} \mathbf{p}_{1}'} = \dots = \overline{\mathbf{o} \mathbf{p}_{n-1}'}$, it holds that
       \begin{equation}
              \overline{\mathbf{p} \mathbf{q}} = n \overline{\mathbf{o} \mathbf{p}_{0}'} + \mathbf{p}_0, \quad (n \in \mathbb{N})
       \end{equation}
       where, $\overline{\mathbf{o} \mathbf{p}_{0}'}$ is prime. As a result, it is sufficient to prove the existence of a decomposition for the case of prime segments. Therefore, we consider the case of prime segments $\overline{\mathbf{o}\mathbf{p}}$. According to Corollary~\ref{CorOfEdgeToShorterAndTriangles}, for any prime segment $\overline{\mathbf{o}\mathbf{p}}$ has a decomposition
       \begin{equation}
              \overline{\mathbf{o}\mathbf{p}} = \triangle \mathbf{o} \mathbf{p} \mathbf{q} + \triangle \mathbf{o} \mathbf{p} \mathbf{r} - \overline{\mathbf{o}\mathbf{q}} - \overline{\mathbf{o}\mathbf{r}},
       \end{equation}
       where $\triangle \mathbf{o} \mathbf{p} \mathbf{q}$ and $\triangle \mathbf{o} \mathbf{p} \mathbf{r}$ are minimum, $\overline{\mathbf{o}\mathbf{q}}$ and $\overline{\mathbf{o}\mathbf{r}}$ are prime and shorter than $\overline{\mathbf{o}\mathbf{p}}$. Since $\overline{\mathbf{o}\mathbf{q}}$ and $\overline{\mathbf{o}\mathbf{r}}$ are prime, we can apply Corollary~\ref{CorOfEdgeToShorterAndTriangles} again. Then the relation
       \begin{equation}
              \begin{aligned}
                     \overline{\mathbf{o}\mathbf{p}} &= \triangle \mathbf{o} \mathbf{p} \mathbf{q} + \triangle \mathbf{o} \mathbf{p} \mathbf{r} - \overline{\mathbf{o}\mathbf{q}} - \overline{\mathbf{o}\mathbf{r}} \\
                     &= \triangle \mathbf{o} \mathbf{p} \mathbf{q} + \triangle \mathbf{o} \mathbf{p} \mathbf{r} - (\mathcal{T}_1' + \mathcal{T}_2' - \overline{\mathbf{oq'}} - \overline{\mathbf{oq''}}) - (\mathcal{T}_3' + \mathcal{T}_4' - \overline{\mathbf{or'} - \mathbf{or''}}) \\
                     &= \sum_{\mathcal{T}_i' \in S'} k_{\mathcal{T}_i'} \mathcal{T}_i' + \overline{\mathbf{oq'}} + \overline{\mathbf{oq''}} + \overline{\mathbf{or'}} + \overline{\mathbf{oq''}},
              \end{aligned}
       \end{equation}
       is obtained where $\mathcal{T}_1', \dots, \mathcal{T}_4', \mathcal{T}_i'$ are minimum triangles, $k_{\mathcal{T}_i'}$ are integers, $S'$ is a finite set of minimum triangles, $\overline{\mathbf{oq'}}$ and $\overline{\mathbf{oq''}}$ are shorter than $\overline{\mathbf{oq}}$, $\overline{\mathbf{or'}}$ and $\overline{\mathbf{or''}}$ are shorter than $\overline{\mathbf{or}}$. Again, we can apply Corollary~\ref{CorOfEdgeToShorterAndTriangles} to $\overline{\mathbf{oq'}}$, $\overline{\mathbf{oq''}}$, $\overline{\mathbf{or'}}$, $\overline{\mathbf{or''}}$, and obtain the form
       \begin{equation}
              \overline{\mathbf{op}} = \sum k_{\mathcal{T}''} \mathcal{T}'' + \sum k_{\mathbf{s}_i} \overline{\mathbf{o}\mathbf{s}_i},
       \end{equation}
       where the number of terms in sums are finite, $k_{\mathcal{T}''}$ and $k_{\mathbf{s}_i}$ are integers, and each $\overline{\mathbf{o}\mathbf{s}_i}$ is prime. This procedure can be repeated as long as the lengths of appearing segments are larger than $1$. Since the lengths of appearing segments strictly decrease at each step, and the minimum length among prime segments is $1$, this iteration must stop within a finite number of steps and the result is expressed as
       \begin{equation}
              \overline{\mathbf{op}} = \sum k_{\mathbf{s}_i'} \overline{\mathbf{o} \mathbf{s}_i'} + \sum k_{\mathcal{T}'''} \mathcal{T}''', \label{EqProofOfSegmentDecomp}
       \end{equation}
       where $k_{\mathbf{s}_i'}$ and $k_{\mathcal{T}'''}$ are integers, and $\overline{\mathbf{o} \mathbf{s}_i'}$ and $\mathcal{T}'''$ denote segments of length $1$ and minimum triangles respectively. Any minimum triangle is obtained by translating a unit triangle, and any segment of length $1$ is obtained by translating $\mathcal{I}_x$ or $\mathcal{I}_y$. Therefore \eqref{EqProofOfSegmentDecomp} can be rewritten as
       \begin{equation}
              \overline{\mathbf{op}} = \mathbf{t} + k_x \mathcal{I}_x + k_y \mathcal{I}_y + \sum_{\mathcal{T} \in S} k_{\mathcal{T}} \mathcal{T},
       \end{equation}
       for a specific translation $\mathbf{t}$.
\end{proof}
We prove Proposition~\ref{PropOfDecompFormula} for the case of polygons as follows. Consider an integral polygon $\mathcal{P}$ and consider a set $P$ which consists of all integral points in $\mathcal{P}$. We can easily obtain a partition $\{P_i\}_{i \in I}$ of $\mathcal{P}$ where $P_i \subset P$ and $\mathrm{conv}(P_i)$ is a minimum triangle for any $i \in I$. Though such a partition is not uniquely determined in general, consider one of the partitions and call it $\Pi_{\min}$. By combining Corollary~\ref{CorOfGeneralSubdivision} and Lemma~\ref{LemOfSegmentDecomposition}, we can prove Proposition~\ref{PropOfDecompFormula}.
\begin{proof}
       There exists at least one minimal partition $\Pi_{\min}$ for any integral polygon $\mathcal{P}$. From this partition, we can obtain a decomposition of $\mathcal{P}$ as
       \begin{equation} \label{Eq1OfProofOfPropOfDecompfomula}
              \mathcal{P} = \sum_{\mathbf{p}_i \in \mathrm{int}(P)} \mathbf{p}_i - \sum_{\overline{\mathbf{p}_i \mathbf{p}_j} \in d(\Pi_{\min})} \overline{\mathbf{p}_i \mathbf{p}_j} + \sum_{P_i \in \Pi_{\min}} \mathrm{conv}(P_i),
       \end{equation}
       by applying Corollary~\ref{CorOfGeneralSubdivision}. Since each $\overline{\mathbf{p}_i \mathbf{p}_j}$ is prime, $\overline{\mathbf{p}_i \mathbf{p}_j}$ can be decomposed as
       \begin{equation} \label{Eq2OfProofOfPropOfDecompfomula}
              \overline{\mathbf{p}_i \mathbf{p}_j} = \mathbf{r}_i + k_{x_i} \mathcal{I}_x + k_{y_i} \mathcal{I}_y + \sum k_{\mathcal{T}} \mathcal{T},
       \end{equation}
       by Lemma~\ref{LemOfSegmentDecomposition}. Substituting \eqref{Eq2OfProofOfPropOfDecompfomula} to \eqref{Eq1OfProofOfPropOfDecompfomula}, we obtain Proposition~\ref{PropOfDecompFormula}.
\end{proof}
\begin{exa}
       We now illustrate the procedure of deriving \eqref{EqOfExamplesOfUnitDecomp}
       \begin{equation}
              \begin{tikzpicture}[baseline = 0.6em, x = 1.4em, y = 1.4em]
                     \coordinate (1) at (0, 0);
                     \coordinate (3) at (1, 1);
                     \coordinate (4) at (3, 0);
                     \fill[black!15] (1) -- (3) -- (4) -- cycle;
                     \draw[step = 1, black, dashed] (-0.2,-0.2) grid (3.2,1.2);
                     \draw (1) -- (3) -- (4) -- cycle;
                     \foreach \P in {1,3,4} \fill[black] (\P) circle (0.11);
              \end{tikzpicture} =
              \begin{tikzpicture}[baseline = 0.6em, x = 1.4em, y = 1.4em]
                     \coordinate (1) at (0, 0);
                     \coordinate (2) at (1, 0);
                     \draw[step = 1, black, dashed] (-0.2,-0.2) grid (1.2,1.2);
                     \draw[thick] (1) -- (2);
                     \foreach \P in {1,2} \fill[black] (\P) circle (0.11);
              \end{tikzpicture} +
              \begin{tikzpicture}[baseline = 0.6em, x = 1.4em, y = 1.4em]
                     \coordinate (1) at (0, 0);
                     \coordinate (2) at (1, 0);
                     \coordinate (3) at (1, 1);
                     \fill[black!15] (1) -- (2) -- (3) -- cycle;
                     \draw[step = 1, black, dashed] (-0.2,-0.2) grid (1.2,1.2);
                     \draw (1) -- (2) -- (3) -- cycle;
                     \foreach \P in {1,2,3} \fill[black] (\P) circle (0.11);
              \end{tikzpicture} +
              \begin{tikzpicture}[baseline = 0.6em, x = 1.4em, y = 1.4em]
                     \coordinate (1) at (1, 0);
                     \coordinate (2) at (0, 1);
                     \coordinate (3) at (2, 0);
                     \fill[black!15] (1) -- (2) -- (3) -- cycle;
                     \draw[step = 1, black, dashed] (-0.2,-0.2) grid (2.2,1.2);
                     \draw (1) -- (2) -- (3) -- cycle;
                     \foreach \P in {1,2,3} \fill[black] (\P) circle (0.11);
              \end{tikzpicture} -
              \begin{tikzpicture}[baseline = 0.6em, x = 1.4em, y = 1.4em]
                     \coordinate (1) at (1, 0);
                     \coordinate (2) at (0, 1);
                     \coordinate (3) at (1, 1);
                     \fill[black!15] (1) -- (2) -- (3) -- cycle;
                     \draw[step = 1, black, dashed] (-0.2,-0.2) grid (1.2,1.2);
                     \draw (1) -- (2) -- (3) -- cycle;
                     \foreach \P in {1,2,3} \fill[black] (\P) circle (0.11);
              \end{tikzpicture}, \nonumber
       \end{equation}
       in Example~\ref{ExamplesOfUnitDecomp}. From the partition relation~\eqref{EqofCorOfGeneralSubdivision}, we can obtain
       \begin{equation}
              \begin{tikzpicture}[baseline = 0.6em, x = 1.4em, y = 1.4em]
                     \coordinate (1) at (0, 0);
                     \coordinate (3) at (1, 1);
                     \coordinate (4) at (3, 0);
                     \coordinate (5) at (1, 0);
                     \coordinate (6) at (2, 0);
                     \fill[black!15] (1) -- (3) -- (4) -- cycle;
                     \draw[step = 1, black, dashed] (-0.2,-0.2) grid (3.2,1.2);
                     \draw (1) -- (3) -- (4) -- cycle;
                     \draw (5) -- (3);
                     \draw (6) -- (3);
                     \foreach \P in {1,3,4,5,6} \fill[black] (\P) circle (0.11);
              \end{tikzpicture} =
              \begin{tikzpicture}[baseline = 0.6em, x = 1.4em, y = 1.4em]
                     \coordinate (1) at (0, 0);
                     \coordinate (2) at (1, 0);
                     \coordinate (3) at (1, 1);
                     \fill[black!15] (1) -- (2) -- (3) -- cycle;
                     \draw[step = 1, black, dashed] (-0.2,-0.2) grid (1.2,1.2);
                     \draw (1) -- (2) -- (3) -- cycle;
                     \foreach \P in {1,2,3} \fill[black] (\P) circle (0.11);
              \end{tikzpicture} +
              \begin{tikzpicture}[baseline = 0.6em, x = 1.4em, y = 1.4em]
                     \coordinate (1) at (1, 0);
                     \coordinate (2) at (0, 1);
                     \coordinate (3) at (0, 0);
                     \fill[black!15] (1) -- (2) -- (3) -- cycle;
                     \draw[step = 1, black, dashed] (-0.2,-0.2) grid (1.2,1.2);
                     \draw (1) -- (2) -- (3) -- cycle;
                     \foreach \P in {1,2,3} \fill[black] (\P) circle (0.11);
              \end{tikzpicture} +
              \begin{tikzpicture}[baseline = 0.6em, x = 1.4em, y = 1.4em]
                     \coordinate (1) at (1, 0);
                     \coordinate (2) at (0, 1);
                     \coordinate (3) at (2, 0);
                     \fill[black!15] (1) -- (2) -- (3) -- cycle;
                     \draw[step = 1, black, dashed] (-0.2,-0.2) grid (2.2,1.2);
                     \draw (1) -- (2) -- (3) -- cycle;
                     \foreach \P in {1,2,3} \fill[black] (\P) circle (0.11);
              \end{tikzpicture} -
              \begin{tikzpicture}[baseline = 0.6em, x = 1.4em, y = 1.4em]
                     \coordinate (1) at (0, 0);
                     \coordinate (2) at (0, 1);
                     \draw[step = 1, black, dashed] (-0.2,-0.2) grid (1.2,1.2);
                     \draw[thick] (1) -- (2);
                     \foreach \P in {1,2} \fill[black] (\P) circle (0.11);
              \end{tikzpicture} -
              \begin{tikzpicture}[baseline = 0.6em, x = 1.4em, y = 1.4em]
                     \coordinate (1) at (0, 1);
                     \coordinate (2) at (1, 0);
                     \draw[step = 1, black, dashed] (-0.2,-0.2) grid (1.2,1.2);
                     \draw[thick] (1) -- (2);
                     \foreach \P in {1,2} \fill[black] (\P) circle (0.11);
              \end{tikzpicture}
       \end{equation}
       with proper translations. According to Corollary~\ref{CorOfEdgeToShorterAndTriangles},
       \begin{equation}
              \begin{tikzpicture}[baseline = 0.6em, x = 1.4em, y = 1.4em]
                     \coordinate (1) at (0, 1);
                     \coordinate (2) at (1, 0);
                     \draw[step = 1, black, dashed] (-0.2,-0.2) grid (1.2,1.2);
                     \draw[thick] (1) -- (2);
                     \foreach \P in {1,2} \fill[black] (\P) circle (0.11);
              \end{tikzpicture} =
              \begin{tikzpicture}[baseline = 0.6em, x = 1.4em, y = 1.4em]
                     \coordinate (1) at (1, 0);
                     \coordinate (2) at (0, 1);
                     \coordinate (3) at (0, 0);
                     \fill[black!15] (1) -- (2) -- (3) -- cycle;
                     \draw[step = 1, black, dashed] (-0.2,-0.2) grid (1.2,1.2);
                     \draw (1) -- (2) -- (3) -- cycle;
                     \foreach \P in {1,2,3} \fill[black] (\P) circle (0.11);
              \end{tikzpicture} +
              \begin{tikzpicture}[baseline = 0.6em, x = 1.4em, y = 1.4em]
                     \coordinate (1) at (1, 0);
                     \coordinate (2) at (0, 1);
                     \coordinate (3) at (1, 1);
                     \fill[black!15] (1) -- (2) -- (3) -- cycle;
                     \draw[step = 1, black, dashed] (-0.2,-0.2) grid (1.2,1.2);
                     \draw (1) -- (2) -- (3) -- cycle;
                     \foreach \P in {1,2,3} \fill[black] (\P) circle (0.11);
              \end{tikzpicture} -
              \begin{tikzpicture}[baseline = 0.6em, x = 1.4em, y = 1.4em]
                     \coordinate (1) at (0, 0);
                     \coordinate (2) at (1, 0);
                     \draw[step = 1, black, dashed] (-0.2,-0.2) grid (1.2,1.2);
                     \draw[thick] (1) -- (2);
                     \foreach \P in {1,2} \fill[black] (\P) circle (0.11);
              \end{tikzpicture} -
              \begin{tikzpicture}[baseline = 0.6em, x = 1.4em, y = 1.4em]
                     \coordinate (1) at (0, 0);
                     \coordinate (2) at (0, 1);
                     \draw[step = 1, black, dashed] (-0.2,-0.2) grid (1.2,1.2);
                     \draw[thick] (1) -- (2);
                     \foreach \P in {1,2} \fill[black] (\P) circle (0.11);
              \end{tikzpicture}
       \end{equation}
       holds. Then we can obtain \eqref{EqOfExamplesOfUnitDecomp} by substitution.
\end{exa}
\section{Max-plus functions and convex polygons} \label{Sec:maxplusfunctions}
A max-plus function, which will be introduced in this section, can be associated with a convex polygon. This connection between max-plus functions and convex polygons enables us to interpret calculations of max-plus functions as geometric operations on convex polygons, such as Minkowski sums. By combining this interpretation and Proposition~\ref{PropOfDecompFormula}, factorizations of a max-plus polynomial with integer coefficients can be obtained.
\subsection{Max-plus function and its calculations}
The maximum of two real numbers $a$ and $b$ is denoted as
\begin{equation*}
       \max(a, b) := 
       \begin{cases}
              a & (a \ge b) \\
              b & (a < b)
       \end{cases}.
\end{equation*}
The following relations hold for any $a, b, c \in \mathbb{R}$,
\begin{gather}
       \max(\max(a, b), c) = \max(a, \max(b, c)), \label{max_asso}\\
       \max(a, b) + c = \max(a + c, b + c) \label{max_dist},
\end{gather}
where $+$ is the standard summation of real numbers. Additionally, it holds that for any $c > 0$
\begin{equation}
       c \max(a,b) = \max(ca, cb) \label{eq:max_Product}.
\end{equation}
We write $\max(a, b, c)$ instead of $\max(\max(a, b), c)$, and analogously for the cases of more than 3 terms. Notation such as $\max_{i \in I}(a_i)$ or $\max_i(a_i)$ are also employed. In relation to the nested expression of two max operations, there are two cases; $\max(\max(a,b),c)$ and $\max(-\max(a,b),c)$. For the former case, we have $\max(a,b,c)$ as above. For the latter case, we derive the following using \eqref{max_asso} and \eqref{max_dist} to expand the nested terms.
\begin{equation}
       \begin{aligned}
              &\max(-\max(a, b), c) \\
              =& \max(-\max(a, b), c) + \max(a, b) - \max(a, b) \\
              =& \max(-\max(a, b) + \max(a, b), c + \max(a, b)) - \max(a, b) \\
              =& \max(0, \max(a + c, b + c)) - \max(a, b) \\
              =& \max(0, a + c, b + c) - \max(a, b).
       \end{aligned} \label{CalcOfNest}
\end{equation}
As a result of these calculations, $\max(-\max(a, b), c)$ is reduced to the difference of two maximum terms $\max(0, a + c, b + c) - \max(a, b)$. It is also possible to rewrite expressions involving multiply nested terms, such as $\max(-\max(-\max(\dots)\dots)\dots)$ or $\max(-\max(\dots),-\max(\dots),\dots)$, as $\max_i(a_i) - \max_j (a_j ')$, where $a_i, a_j' \in \mathbb{R}$, via similar calculations to \eqref{CalcOfNest}.

We investigate piecewise-linear functions in the form of $f(x, y) = \max_i(a_i x + b_i y): \mathbb{R}^2 \to \mathbb{R}$, where $x$ and $y$ are real variables, and the coefficients $a_i, b_i$ are real constants. Let us denote the class of such functions by $\mathrm{MP}$. We can verify that $\mathrm{MP}$ is closed under $\max$ and $+$. By \eqref{max_asso} and \eqref{max_dist}, it holds that
\begin{align}
       \max(\max_{i \in I}(a_i x + b_i y), \max_{j \in J}(a_j' x + b_j' y)) = &\max_{k \in I \cup J}(a_k'' x + b_k'' y), \label{EquiOfMPMax}\\
       \max_{i \in I}(a_i x + b_i y) + \max_{j \in J}(a_j' x + b_j' y) = &\max_{(i,j) \in I \times J}((a_i + a_j')x + (b_i + b_j')y). \label{EquiOfMPPlus}
\end{align}
Applying \eqref{EquiOfMPMax} and \eqref{EquiOfMPPlus}, various functions consisting of multiple $\max$ and $+$ can be reduced to the form $\max_i (a_ix + b_i y)$ and it turns out that such functions belong to $\mathrm{MP}$.

$\mathrm{MP}$ is not closed under the operation $-$. For instance, functions such as $-f$, $\max(f, -g)$, and $-\max(f, g)$ may not be included in $\mathrm{MP}$ even if $f, g \in \mathrm{MP}$. In general, $\max(\max_i(a_i x + b_i y), -\max_j(a_j' x + b_j' y))$ cannot be written in the form of a function $\max_k (a_k x + b_k y)$. However, it can be reduced to the difference of functions of $\mathrm{MP}$ by similar calculations to \eqref{CalcOfNest} via
\begin{equation}
       \begin{aligned}
              &\max(\max_i(a_i x + b_i y), -\max_j(a_j' x + b_j' y)) \\
              = &\max(\max_i(a_i x + b_i y) + \max_j(a_j' x + b_j' y), 0) - \max_j(a_j' x + b_j' y) \\
              = &\max(\max_{(i,j) \in I \times J}((a_i + a_j')x + (b_i + b_j')y), 0x + 0y) - \max_j(a_j' x + b_j' y).
       \end{aligned}
\end{equation}
Similarly, functions in more complicated expressions involving $\max$, $+$, $-$, and $a_k x + b_k y$ can be rewritten as $\Sigma \max_i (a_i x + b_i y) - \Sigma \max_j (a_j' x + b_j' y)$ with specific $a_i, b_i, a_j', b_j'$, where $\Sigma \max_i (a_i x + b_i y) \in \mathrm{MP}$ and $\Sigma \max_j (a_j' x + b_j' y) \in \mathrm{MP}$.

Two functions of $\mathrm{MP}$ with different coefficients $\{(a_i, b_i) \,;\, i \in I\}$ and $\{(a_j', b_j') \,;\, j \in J\}$ may have the same value at any $(x,y) \in \mathbb{R}^2$ and be equivalent to each other as shown in the following example.
\begin{exa} \label{ExaOfInnerPoint}
       Assume $f_1(x,y) := \max(0, 2x, y)$ and $f_2(x,y) := \max(0, x, 2x, y)$. The term $x$ in $f_2$ is smaller than $2x$ if $x > 0$, and otherwise smaller than or equal to $0$. Consequently, the term $x$ does not affect the value of $f_2$. Therefore, $\max(0, 2x, y) = \max(0, x, 2x, y)$ at any $(x,y) \in \mathbb{R}^2$ and $f_1$ is equivalent to $f_2$. Similarly, we can easily see that infinitely many functions are equivalent to $f_1$, such as $\max(0,2x, 2c_1x, y)$ where $(0 < c_1 < 1)$, $\max(0,2x, 2c_2x + (1 - c_2)y, y)$ where $(0 < c_2 < 1)$, and so on. The function $f_1$ includes $3$ terms in a maximum operation and we note that it is the least number of terms necessary to express the function equivalent to $f_1$ in the form of $\max_i(a_i x + b_i y)$.
\end{exa}
\subsection{Interpretations of max-plus functions as convex polygons}
For a finite set $P = \{\mathbf{p}_1 = (a_1, b_1), \dots, \mathbf{p}_m = (a_m, b_m)\}$ of points, we employ $\max \langle P\rangle$ to denote a function of $\mathrm{MP}$ such that $\max \langle P\rangle := \max(a_1 x + b_1 y, \dots, a_m x + b_m y)$. The following lemma then holds.
\begin{lem} \label{LemOfVer}
       Two functions $\max \langle P\rangle$ and $\max \langle v(P)\rangle$ are equivalent in that
       \begin{equation}
              \max(a_1 x + b_1 y, \dots, a_m x + b_m y) = \max(a_{v_1} x + b_{v_1} y, \dots, a_{v_l} x + b_{v_l} y)
       \end{equation}
       for any $(x,y) \in \mathbb{R}^2$.
\end{lem}
Note that $v(P)$ is defined in Section~\ref{Sec:PolytopesAndConvexGeometry} as a set of vertices of $\mathrm{conv}(P)$.
\begin{proof}
       Suppose $(a_i,b_i) \in P$ and $(a_i,b_i) \notin v(P)$. Fix a point $(x',y') \in \mathbb{R}^2$. Then $a_i x' + b_i y'$ can be written as
       \begin{equation}
              a_i x' + b_i y' = \Sigma_{j = 1}^{l} c_j (a_{v_j} x' + b_{v_j} y'),
       \end{equation}
       where each $c_j \ge 0$ and $\Sigma_{j = 1}^{l} c_j = 1$, according to \eqref{EquiOfVertex}. There exists a point $(a_{v_M},b_{v_M}) \in v(P)$ satisfying
       \begin{equation}
              \max_{(a_{v_i},b_{v_i}) \in v(P)}(a_{v_1} x' + b_{v_1} y', \dots, a_{v_l} x' + b_{v_l} y') = a_{v_M} x' + b_{v_M} y'.
       \end{equation}
       Then, we obtain
       \begin{equation}
              \begin{aligned}
                     a_i x' + b_i y' &= \Sigma_{j = 1}^{l} c_j (a_{v_j} x' + b_{v_j} y') \\
                     &\le \Sigma_{j = 1}^{l} c_j (a_{v_M} x' + b_{v_M} y') \\
                     &= a_{v_M} x' + b_{v_M} y'.
              \end{aligned}
       \end{equation}
       Similarly, we can prove that at least one of $a_{v_j} x + b_{v_j} y$ is greater than $a_i x + b_i y$ for any $(x,y) \in \mathbb{R}^2$. Therefore, we can remove the term $a_ix+b_iy$ from $\max \langle P\rangle$. Removing all terms $a_ix+b_iy$ such that $(a_i,b_i) \notin v(P)$ from $\max \langle P\rangle$, we obtain $\max \langle v(P)\rangle$. Thus, $\max \langle v(P)\rangle$ is equivalent to $\max \langle P\rangle$.
\end{proof}
According to Lemma~\ref{LemOfVer}, we can judge whether or not two functions $\max_i(a_i x + b_i y)$ and $\max_j(a_j' x + b_j' y)$ are equivalent by comparing the convex hull of points $(a_i, b_i)$ and $(a_j', b_j')$.
\begin{cor} \label{CorOfEqui}
Let $P_1$ and $P_2$ be finite sets of $\mathbb{R}^2$, then $\max\langle P_1\rangle$ and $\max\langle P_2\rangle$ are equivalent if and only if $\mathrm{conv}(P_1) = \mathrm{conv}(P_2)$. \label{equi1}
\end{cor}
\begin{proof}
       It immediately follows from Lemma~\ref{LemOfVer}.
\end{proof}
We now give an example of Corollary~\ref{equi1}.
\begin{exa}
       Let $P_1 = \{(0,0),(2,0),(0,1)\}$ and $P_2 = \{(0,0),(1,0),(2,0),(0,1)\}$. Both of $\mathrm{conv}(P_1)$ and $\mathrm{conv}(P_2)$ are the same triangle $\mathrm{conv}\{(0,0),(2,0),(0,1)\}$. On the other hand, as already seen in Example~\ref{ExaOfInnerPoint},
       \begin{equation}
              \max \langle P_1\rangle = \max(0,2x,y) = \max(0,x,2x,y) = \max \langle P_2\rangle, \label{ExaOfEqu}
       \end{equation}
       for any $(x,y) \in \mathbb{R}^2$.
\end{exa}
Combining Corollary~\ref{equi1} and the geometric identities \eqref{rel2} and \eqref{rel3}, we can derive the following identities.
\begin{cor} \label{CorOfCorMaxPol}
       For any finite sets of points $P_1$ and $P_2$ of $\mathbb{R}^2$, it holds that
       \begin{gather}
              \max \langle P_1 \cup P_2\rangle = \max(\max \langle P_1\rangle, \max \langle P_2\rangle), \label{equi2}\\
              \max \langle P_1 + P_2\rangle = \max \langle P_1\rangle + \max \langle P_2\rangle. \label{equi3}
       \end{gather}
       Additionally, if the Minkowski difference $\mathrm{conv}(P_1) - \mathrm{conv}(P_2) =: \mathrm{conv}(P_3)$ exists, then
       \begin{equation}
              \max \langle P_1 \rangle - \max \langle P_2 \rangle = \max \langle P_3 \rangle, \label{Eq3OfCorOfCorMaxPol}
       \end{equation}
       holds.
\end{cor}
\begin{proof}
       Equations \eqref{equi2} and \eqref{equi3} are immediately derived from the identities \eqref{rel2}, \eqref{rel3} and Corollary~\ref{CorOfEqui}. For \eqref{Eq3OfCorOfCorMaxPol}, we have $\mathrm{conv}(P_1) = \mathrm{conv}(P_2) + \mathrm{conv}(P_3)$ according to the definition of the Minkowski difference described in Section~\ref{Sec:PolytopesAndConvexGeometry}. Since $\mathrm{conv} (P_2) + \mathrm{conv}(P_3) = \mathrm{conv}(P_2 + P_3)$, we have $\mathrm{conv}(P_1) = \mathrm{conv}(P_2 + P_3)$. Then $\max \langle P_1\rangle = \max \langle P_1 + P_2\rangle = \max \langle P_1\rangle + \max \langle P_2\rangle$ is obtained by Corollary~\ref{CorOfEqui} and \eqref{equi3}.
\end{proof}
We now provide an example of Corollary~\ref{CorOfCorMaxPol}.
\begin{exa}
Assume $P_1 = \{(0,0), (1,1), (1,2)\}$ and $P_2 = \{(0,0), (1,1), (2,0)\}$. The union of $P_1$ and $P_2$ is $P_1 \cup P_2 = \{(0,0), (1,1), (2,0), (1,2)\}$, and thus $\max \langle P_1 \cup P_2\rangle$ is
\begin{equation}
       \max \langle P_1 \cup P_2\rangle = \max(0, x+y, 2x, x + 2y).
\end{equation}
On the other hand, according to \eqref{max_asso}, the function $\max(\max \langle P_1\rangle, \max\langle P_2\rangle)$ can be rewritten as
\begin{equation}
\begin{aligned}
\max(\max \langle P_1\rangle, \max\langle P_2\rangle) &= \max(\max(0, x + y, x + 2y), \max(0, x + y, 2x)) \\
&= \max(0, x+y, 2x, x + 2y).
\end{aligned} \label{ExaOfMax}
\end{equation}
Thus, we obtain that $\max(\max \langle P_1\rangle, \max\langle P_2\rangle)$ and $\max\langle P_1 \cup P_2\rangle$ are equivalent. Similarly, assume $P_1 = \{(0,0),(1,0),(1,1)\}$ and $P_2 = \{(0,0),(0,1),(1,1)\}$, then 
\begin{equation}
       P_1 + P_2 = \{(0,0),(1,0),(0,1),(1,1),(2,1),(1,2),(2,2)\}.
\end{equation}
On the other hand, $\max\langle P_1 \rangle + \max\langle P_2 \rangle$ can be calculated as follows,
\begin{equation}
\begin{aligned}
&\max\langle P_1\rangle + \max\langle P_2\rangle \\
=& \max(0,x,x+y) + \max(0,y,x+y) \\
=& \max(\max(0,x,x+y),\max(y,x+y,x+2y),\max(x+y,2x+y,2x+2y)) \\
=& \max(0, x, y, x + y, 2x + y, x + 2y, 2x + 2y),
\end{aligned} \label{ExaOfPlus}
\end{equation}
and the last of \eqref{ExaOfPlus} corresponds to $\max\langle P_1 + P_2\rangle$.
\end{exa}
Through Proposition~\ref{LemOfVer}, Corollary~\ref{CorOfEqui}, and Corollary~\ref{CorOfCorMaxPol}, we have seen that functions of $\mathrm{MP}$ obey formulae analogous to those obtained for sets of points on $\mathbb{R}^2$ and convex hulls of those sets. Each function of $\mathrm{MP}$ can be associated with a convex polygon, which is a convex hull of points whose coordinates are equal to coefficients of terms in the function. Two operations $\max$ and $+$ acting on functions of $\mathrm{MP}$ can be calculated by the union and the Minkowski sum of sets of points respectively. These geometric interpretations reveal the essence of max-plus functions. Indeed, if an identity involving Minkowski sums or differences of convex polygons, such as the cutting relation and the partition relation, is given, it induces an identity involving sums or differences of max-plus functions associated with the polygons. Considering these correspondences between $\mathrm{MP}$ functions and convex polygons, we can decompose a given $\mathrm{MP}$ function into a sum of ``prime'' $\mathrm{MP}$ functions using Proposition~\ref{PropOfDecompFormula}.
\begin{cor}
       For any finite set $P \subset \mathbb{Z}^2$, $\max\langle P \rangle$ can be expressed as
       \begin{equation}
              \max\langle P \rangle = a_0 x + b_0 y + k_x \max(0,x) + k_y \max(0,y) + \sum_{\mathcal{T} \in S} k_{\mathcal{T}} \max\langle v(\mathcal{T})\rangle, \label{EqOfCorMPPolFactor}
       \end{equation}
       where $(a_0,b_0) \in \mathbb{Z}^2$, $k_x$, $k_y$, $k_{\mathcal{T}}$ are integers, and $S$ is a finite subset of $T_{\mathrm{uni}}$. \label{CorMPPolFactor}
\end{cor}
\begin{proof}
       The polygon $\mathrm{conv}(P)$ is an integral polygon. Therefore, we can obtain a decomposition of $\mathrm{conv}(P)$ by applying Proposition~\ref{PropOfDecompFormula} such that
       \begin{equation}
              \mathrm{conv}(P) = (a_0,b_0) + k_x \mathcal{I}_x + k_y \mathcal{I}_y + \sum_{\mathcal{T} \in S} k_{\mathcal{T}} \mathcal{T}, \label{EqOfProofOfCorMPPolFactor}
       \end{equation}
       where $(a_0, b_0) \in \mathbb{Z}^2$, $k_x$, $k_y$, $k_\mathcal{T}$ are integers, and $S$ is a finite subset of $T_{\mathrm{uni}}$. By applying Corollary~\ref{CorOfEqui} to \eqref{EqOfProofOfCorMPPolFactor}, we obtain \eqref{EqOfCorMPPolFactor}.
\end{proof}
\begin{rem}
       By regarding $\max$ as addition $\oplus$ and $+$ as multiplication $\otimes$, respectively, the expression $\max_i(a_i x + b_i y + c_i)$ can be rewritten as a tropical polynomial in the form of ``$\bigoplus_i c_i \otimes x^{\otimes a_i} \otimes y^{\otimes b_i}$''. The function $\max \langle P \rangle$ appearing in Corollary~\ref{CorMPPolFactor} can be rewritten as a tropical polynomial ``$\bigoplus_i x^{\otimes a_i} y^{\otimes b_i}$'' such that all coefficients of monomials are $0$, which is the unity of the tropical semifield. As a consequence of Corollary~\ref{CorMPPolFactor}, it is revealed that the set of all tropical polynomials with coefficients $c_i = 0$ are generated by $x$, $y$, $\max(0,x)$, $\max(0,y)$, and $\{\max\langle v(\mathcal{T})\rangle \,|\, \mathcal{T} \in T_{\mathrm{uni}}\}$. Note also that the polygon $\mathrm{conv}\{(a_i, b_i)\}_i$ corresponds to the Newton polygon of a tropical polynomial ``$\bigoplus_i c_i \otimes x^{\otimes a_i} \otimes y^{\otimes b_i}$''. Newton polygons are regarded as fundamental objects in the study of algebraic properties of tropical polynomials.
\end{rem}
According to Corollary~\ref{CorMPPolFactor}, a difference $\max \langle P' \rangle - \max \langle P'' \rangle$, where $P', P'' \subset \mathbb{Z}^2$, can also be expressed in terms of $a_0 x + b_0 y$, $\max(0,x)$, $\max(0,y)$, and $\max\langle v(\mathcal{T})\rangle$.
\begin{cor}
       For any pair of finite sets $P', P'' \subset \mathbb{Z}^2$, $\max \langle P' \rangle - \max \langle P'' \rangle$ can be expressed as
       \begin{equation}
              \max \langle P' \rangle - \max \langle P'' \rangle = a_0 x + b_0 y + k_x \max(0,x) + k_y \max(0,y) + \sum_{\mathcal{T} \in S} k_{\mathcal{T}} \max\langle v(\mathcal{T})\rangle, \label{EqOfCorMPRatFactor}
       \end{equation}
       where $(a_0,b_0) \in \mathbb{Z}^2$, $k_x$, $k_y$, $k_{\mathcal{T}}$ are integers, and $S$ is a finite subset of $T_{\mathrm{uni}}$. \label{CorMPRatFactor}
\end{cor}
\begin{proof}
       As a result of Corollary~\ref{CorMPPolFactor}, $\max \langle P' \rangle$ and $\max \langle P'' \rangle$ can be represented as
       \begin{align}
              \max \langle P' \rangle &= a_0' x + b_0' y + k_x' \max(0,x) + k_y' \max(0,y) + \sum_{\mathcal{T}' \in S'} k_{\mathcal{T}'}' \max\langle v(\mathcal{T}')\rangle, \\
              \max \langle P'' \rangle &= a_0'' x + b_0'' y + k_x'' \max(0,x) + k_y'' \max(0,y) + \sum_{\mathcal{T}'' \in S''} k_{\mathcal{T}''}'' \max\langle v(\mathcal{T}'')\rangle,
       \end{align}
       where $S'$ and $S''$ are finite subsets of $T_{\mathrm{uni}}$. Therefore, by taking $a_0 = a_0' - a_0''$, $b_0 = b_0' - b_0''$, and so on, we obtain \eqref{EqOfCorMPRatFactor}.
\end{proof}
\begin{rem}
       Corollary~\ref{CorMPRatFactor} is applicable in the case where the Minkowski difference $\mathrm{conv}(P') - \mathrm{conv}(P'')$ does not exist. Although the expression $\mathcal{P}' - \mathcal{P}''$ cannot be defined for an arbitrary pair of polygons $\mathcal{P}'$ and $\mathcal{P}''$ in the sense of Minkowski differences, a formal treatment of such a difference has been constructed. The formal treatment is referred to as a virtual polytope~\cite{PK1993,PS2015}, and regarded as a generalization of convex polytopes. Originally, cutting relations and partition relations were proven for virtual polytopes. Therefore, Proposition~\ref{PropOfDecompFormula} can be extended to virtual polytopes in a similar way to Corollary~\ref{CorMPRatFactor}. This implies that our results can be extended to virtual polytopes, and we will consider this problem in future work.
\end{rem}
\section{Concluding remarks and discussions} \label{Sec:Concluding}
We have proved that any integral polygon can be expressed as a signed Minkowski decomposition composed of unit triangles and unit segments. The procedure to obtain such a decomposition is as follows. Firstly, we choose a partition of the given polygon with minimum triangles. By applying the partition relations (Corollary~\ref{CorOfGeneralSubdivision}) to this partition, a signed Minkowski decomposition which consists of unit triangles and integral line segments can be obtained. The lengths of line segments appearing in this decomposition are in general longer than $1$. According to Lemma~\ref{LemOfSegmentDecomposition}, any integral line segments can be expressed as a signed decomposition consisting of unit triangles and unit segments. Therefore, the decomposition can be expressed by only using unit triangles and unit segments with a translation. 

Lemma~\ref{LemOfSegmentDecomposition} enables us to obtain finer decompositions than those which have been previously introduced~\cite{Funke2021}, with respect to the size of components. The procedure in~\cite{Funke2021} gives a decomposition of an integral polygon consisting of right-angled triangles whose areas are larger than $1/2$ and integral line segments whose lengths in general are longer than $1$. In other words, we have constructed a generating set of integral polygons with respect to a signed Minkowski decomposition, which consists of elements smaller than those already known~\cite{Funke2021}. This is a noteworthy outcome of our work. Note that Funke's decomposition is applicable for not only $2$-dimensional polygons, but also for $d$-dimensional polytopes $(d \ge 3)$, although our present method is only for $2$-dimensional polygons. Moreover, Funke's generating set forms a basis of integral polytopes, and thus Funke's decomposition is uniquely determined.

According to the correspondence between convex polygons and max-plus functions, we can decompose a max-plus function with integer coefficients into a sum of `prime' max-plus functions of $2$ or $3$ terms. There already exist some results about standard forms of max-plus functions with real coefficients~\cite{BGH2025,KMPS2023,KS1987,TW2024}. Our results demonstrate decompositions of functions with integer coefficients. This is a particular novelty of our work. The procedure to obtain the expression is obvious from the proof of Proposition~\ref{PropOfDecompFormula}, and thus our result may have various applications.
\subsection{Discussions and further works}
We propose the following topics as future work.

\textbf{Uniqueness of decompositions}
Although the uniqueness of decompositions is not verified in the proof of Proposition~\ref{PropOfDecompFormula}, we observed that any $2$-dimensional integral polygons seem to be decomposed uniquely following our procedure. However, the proof of uniqueness is not completed. If the proof is completed, we can obtain another basis of the $2$-dimensional subgroup of the integral polytope group, distinguished from the one discovered by Funke~\cite{Funke2021}.

\textbf{Higher dimensional polytopes}
Another subject is to extend our results to the decompositions for $d$-dimensional polytopes such that $d \ge 3$, but there are some obstacles. First, our procedure to obtain decompositions is based on the fact that every integral polygon have a partition consisting only of triangles of area 1/2. Such partitions are often referred to as unimodular triangulations~\cite{HPPS2021}. For dimension $d \ge 3$, it is known that a unimodular triangulation by $d$-simplices does not exist for any $d$-dimensional convex polytope. Additionally, we have not obtained a higher dimensional version of Lemma~\ref{LemOfSegmentDecomposition}, which remains an open problem.

\textbf{Analysis on complexity and efficiency}
There are a number of previous studies on the complexity and efficiency of Minkowski decompositions of polygons from the viewpoint of computational geometry and engineering~\cite{Fukuda2004,Hachenberger2009}. Furthermore, some studies have emerged aiming at constructing efficient algorithms for factorizing tropical polynomials by utilizing the correspondence between Minkowski decompositions of polygons and factorizations of tropical polynomials~\cite{KMPS2023,MRZ2022,TW2024}. In the future, we may analyze the computational efficiency of our algorithm.
\section*{Acknowledgements}
This work was supported by JST SPRING, Grant Number JPMJSP2128.
\section*{Statements and Declarations}
The author has no conflict of interests to declare that are relevant to the content of this article.

\end{document}